\documentclass[11pt]{amsart}
\usepackage{amsmath}
\usepackage{amscd}
\usepackage{amsthm}
\usepackage{amssymb} \usepackage{latexsym}
\usepackage{eufrak}
\usepackage{euscript}
\usepackage{epsfig}
\usepackage{graphics}
\usepackage{array}
\usepackage{enumerate}
\usepackage{dsfont}
\usepackage{color}

\newcommand{\notequiv}{\equiv \hskip -3.4mm \slash}
\newcommand{\smallnotequiv}{\equiv \hskip -1.7mm \slash}

\newcommand{\bel}[1]{\begin{equation}\label{#1}}

\newcommand{\be}{\begin{equation}}

\newcommand{\ba}{\begin{eqnarray}}
\newcommand{\ea}{\end{eqnarray}}
\newcommand{\rf}[1]{(\ref{#1})}
\newcommand{\bi}{\bibitem}
\newcommand{\qe}{\end{equation}}
\newcommand{\R}{{\mathbb R}}
\newcommand{\N}{{\mathbb N}}

\theoremstyle{theorem}
\newtheorem{theo}{Theorem}[section]
\newtheorem{satz}{Proposition}[section]
\newtheorem{example}{Example}[section]
\theoremstyle{corollary}
\newtheorem{coro}{Corollary}[section]
\theoremstyle{lemma}
\newtheorem{lemma}{Lemma}[section]
\theoremstyle{definition}
\newtheorem{defi}{Definition}[section]
\theoremstyle{proof}

\theoremstyle{remark}
\newtheorem*{rem}{Remark}

\begin{document}

\title[Bipartite and neighborhood graphs and the graph Laplacian]{Bipartite and neighborhood graphs and the spectrum of the
  normalized graph Laplace operator}
\author{Frank Bauer and J\"urgen Jost}
\thanks{Max Planck Institute for Mathematics in the Sciences,
Inselstrasse 22, D-04103 Leipzig, Germany}

\maketitle \textbf{To appear in Communications in Analysis and
Geometry}
\begin{abstract}
We study the spectrum of the normalized Laplace operator of a
connected graph $\Gamma$. As is well known, the smallest
nontrivial eigenvalue measures how difficult it is to decompose
$\Gamma$ into two large pieces, whereas the largest eigenvalue
controls how close $\Gamma$ is to being bipartite. The smallest
eigenvalue can be controlled by the Cheeger constant, and we
establish a dual construction that controls the largest
eigenvalue. Moreover, we find that the neighborhood graphs
$\Gamma[l]$ of order $l\geq2$ encode important spectral
information about $\Gamma$ itself which we systematically explore.
In particular, the neighborhood graph method leads to new
estimates for the smallest nontrivial eigenvalue that can  improve
the Cheeger inequality, as well as an explicit estimate for the
largest eigenvalue from above and below. As applications of such
spectral estimates, we provide a criterion for the
synchronizability of coupled map lattices, and an estimate for the
convergence rate of random walks on graphs.
\end{abstract}
 \tableofcontents
\section{Introduction}
A general  principle in geometry tells us that the spectrum of a
Laplace operator encodes important geometric information about the
underlying space. This principle has been particularly fertile in
Riemannian geometry. One of the key questions has been the control
from below of the first nonzero eigenvalue of the Laplace-Beltrami
operator in terms of the geometry of the underlying Riemannian
manifold (assumed to be compact here for simplicity of
exposition). The Lichnerowicz bound estimates the first eigenvalue
from below in terms of a lower bound for the Ricci curvature. In
contrast, the Cheeger estimate controls the first eigenvalue from
below in terms of a global quantity that expresses how difficult
it is to cut the manifold into two large pieces  \cite{Cheeger70}.
In this way, the first eigenvalue could be related to the
fundamental analytic constants of a Riemannian manifold, like the
isoperimetric or Sobolev constants. The  work of Li and Yau
\cite{Li-Yau80} utilized gradient bounds for eigenfunctions in
order to control the first eigenvalue from below in terms of the
diameter and Ricci bounds of the Riemannian manifold. More
generally, their famous Harnack inequality for the heat kernel
\cite{Li-Yau86} then allowed for a systematic control of all
eigenvalues of a Riemannian manifold, with the optimal asymptotics
as given by Weyl's law. See for instance \cite{Chavel84} for a
systematic treatment of
eigenvalues in Riemannian geometry.\\
In graph theory, the algebraic graph Laplace operator has been
explored for a long time, see for example \cite{Merris94}. More
recently, F.Chung and S.T.Yau, see e.g.
\cite{Chung-Yau95,Chung-Grigoryan-Yau96,Chung-Yau98} and the
monograph \cite{Chung97},  systematically investigated the
normalized graph Laplace operator $\Delta$ of an unweighted and
undirected graph. This operator, which is different from the
algebraic graph Laplace operator, underlies random walks and
diffusion processes with conservation laws on graphs. The
normalized graph Laplace operator is related to the
Laplace-Beltrami operator for a Riemannian manifold. Thus, in
order to study the spectrum of $\Delta$, one can systematically
apply methods developed in Riemannian geometry for the
investigation of the spectrum of the Laplace-Beltrami operator and
this led to many
remarkable insights, see the works just cited and the references therein.\\
In particular,
  the  smallest nontrivial eigenvalue can be well controlled  in terms of the Cheeger
constant \cite{Chung97}. (In a graph-theoretical setting, such
constants can already be found
in  earlier work by Poly\'a and Szeg\"o \cite{Polya51}.)\\
In contrast to a Riemannian manifold, on a graph, the spectrum of
the normalized Laplace operator is always bounded from above. In
fact, the upper bound 2 is achieved if and only if the graph is
bipartite. (We recall that a graph is bipartite if its vertex set
consists of two classes such that edges are only permitted between
two vertices from opposite classes.) Therefore, it is a natural
question how to control the largest eigenvalue for graphs that are
not bipartite. The original goal of  this article was to derive
bounds for the largest eigenvalue of $\Delta$ from above and
below. These bounds reflect how different the graph in question is
from a bipartite one resp. how close it is
to such a graph. \\
In fact, however, these estimates led us to discover more general
structures that go beyond our original goal. First of all, we
construct a dual version $\bar{h}$ of the Cheeger constant $h$ and
derive bounds for the largest eigenvalue from above and below in
terms of $\overline{h}$. We find interesting relations between $h$
and $\bar{h}$, and the combination of these two constants tells us
more about the graph than either of them does individually.
Moreover, we find that the neighborhood graphs $\Gamma[l]$, of
order $l\geq2$, of a graph $\Gamma$ also encode important spectral
information about $\Gamma$ itself. For concreteness let $l=2$ for
the moment.  The idea then is that $\Gamma[2]$ is a weighted graph
with the same vertices as $\Gamma$ itself, and two vertices are
connected in $\Gamma[2]$ when they share at least one neighbor in
$\Gamma$, with lower weights for more shared neighbors. More
precisely, let $\Gamma$ be a weighted graph, with the weight of
the edge between the vertices $i$ and $j$ denoted by $w_{ij}$
(which is 0 unless $i$ and $j$ are neighbors), and the degree of
$i$ being $d_i=\sum_j w_{ij}$. For the neighborhood graph
$\Gamma[2]$, the weight of the edge $e[2]=(i,j)$ in $\Gamma[2]$
then is given by $w_{ij}[2] = \sum_k \frac{1}{d_k}w_{ik}w_{kj}$.
Consequently, $i$ and $j$ are neighbors in $\Gamma[2]$ if $i$ and
$j$ have at least one common neighbor in $\Gamma$, i.e. there
exists a path of length $2$ between $i$ and $j$ in $\Gamma$. Note
that the weights of $\Gamma[2]$ are normalized in such a manner
that every vertex $i$ has the same degree in both $\Gamma$ and
$\Gamma[2]$. It turns out that the eigenvalues of $\Gamma[2]$ are
given by $\lambda (2-\lambda)$ when $\lambda$ stands for the
eigenvalues of $\Gamma$. This well suits our purpose because
controlling the highest eigenvalue $\lambda_{max}$ from above is
equivalent to controlling $2-\lambda_{max}$ from below. Thus,
lower spectral bounds on $\Gamma[2]$ yield lower and upper
spectral bounds for $\Gamma$. For certain graphs, these new lower
bounds improve the Cheeger estimate for the smallest nonzero
eigenvalue.  We not only utilize this principle to derive such
bounds, but we also explore the relation between the spectra of
$\Gamma$ and $\Gamma[2]$ in more general terms. Naturally, the
construction of the neighborhood graph $\Gamma[2]$ can be
generalized to higher order neighborhood graphs, i.e. $i$ and $j$
are neighbors in $\Gamma[l]$ if there exists a path of length $l$
between $i$ and $j$ in $\Gamma$. Again, the weights of $\Gamma[l]$
are normalized in such a way that every vertex has the same degree
in both $\Gamma$ and $\Gamma[l]$. In the present paper, we also
explore the spectra of the higher order neighborhood graphs
$\Gamma[l]$ and their relations to the spectrum of $\Gamma$. The
concept of the neighborhood graph is quite general, and can also
be used to investigate the spectrum of the normalized graph
Laplace operator defined on directed \cite{FB}  graphs.
\\In the last two sections, we will apply our new eigenvalue bounds
to two concrete problems - the convergence of random walks on
graphs, and the synchronization for coupled map lattices, that is,
a dynamical system supported on the vertices of a graph and
coupled according to the interaction structure given by the edges
of the graph. Again, the principle is  that eigenvalue estimates
control how different the graph in question is from the two
extremes of a disconnected or a bipartite graph. On a disconnected
or a bipartite graph, for different reasons, the random walk does
not converge to a stationary distribution, and the coupled map
lattice does not synchronize.

\section{The graph Laplace operator and its basic properties}
In this paper, $\Gamma$ is an undirected, weighted, connected,
finite graph of $N$ vertices. We do not exclude loops, i.e., edges
connecting a vertex with itself. The  vertices are denoted by
$i,j,\dots $. $V$ denotes the vertex  and $E$ the edge set of
$\Gamma$, respectively. When the vertices $i$ and $j$ are
connected by an edge, they are called neighbors, in symbols $i\sim
j$. The associated weight function $w:V\times V \rightarrow
\mathds{R}$ satisfies $w_{ij} = w_{ji}$ and $w_{ij}> 0 $ whenever
$i\sim j$ and $w_{ij} = 0$ iff $i\nsim j$. For a vertex $i$, its
degree $d_i$ is given by $d_i:=\sum_jw_{ij}$. When  $w_{ij}= 1 $
whenever $i\sim j$, we shall speak of an unweighted graph.

The clustering coefficient $C$ of an unweighted graph $\Gamma$ is
defined as \bel{2a.1} C:=\frac{3 \times \mbox{ number of triangles
}}{\mbox{number of
    connected triples of vertices}},
\end{equation}
where  a triangle is a triple of mutually connected vertices. The
clustering coefficient measures how many connections there exist
between the neighbors of a node. $C$ becomes maximal if $\Gamma$
is a fully connected graph. In contrast, $C$ vanishes when
$\Gamma$ is a bipartite graph, that is, consists of two classes
$V_1,V_2$ of vertices such that no vertices in the same class are
connected by an edge. In particular, there are no loops in a
bipartite graph. Equivalently, a graph is bipartite iff it has no
cycles of odd length, and thus in particular no triangles.

We now recall the definition of the normalized graph Laplace
operator and state its basic properties.
\begin{defi}We have a natural measure $\mu$ on the vertex set $V$ given by $\mu(i)=d_i$.
The inner product of two functions $u,v\in \ell^2(V,\mu)$ is
defined  as
 \be
\label{77} (u,v)_\mu := \sum_{i\in V}\mu(i)u(i)v(i).\qe The
Hilbertspace $\ell^2(V,\mu)$ is then given by
\[\ell^2(V,\mu)= \{u:V\rightarrow \mathds{R}\quad| \quad(u,u)_\mu <\infty\}.\]
\end{defi} Since we consider only finite graphs here, the space
$\ell^2(V,\mu)$ is nothing but the space of real-valued functions
on $V$ endowed with the inner product $(\cdot,\cdot)_\mu$. We
study the normalized graph Laplace operator \be \nonumber
\Delta:\ell^2(V,\mu) \rightarrow \ell^2(V,\mu)
\end{equation}
\bel{3} \Delta v(i):= \frac{1}{d_i} \left(\sum_{j}w_{ij}(v(i) -
 v(j))\right).
\end{equation}
The Laplace operator $\Delta$ underlies random walks on graphs. In
fact, the Laplace operator $\Delta$ can be considered as $\Delta
=: I -P$, where $I$ denotes the identity and $P$ is transition
probability operator of a random walk (or sometimes called the
Markov operator), respectively. We should point our here that the
normalized graph Laplace operator $\Delta$ is not exactly the one
studied by Fan Chung \cite{Chung97}. However, both Laplace
operators are unitarily equivalent and therefore  have the same
spectrum. We recall the following basic properties:
\begin{enumerate}
\item[(i)] $\Delta$ is selfadjoint w.r.t. $(.,.)_\mu$, i.e.
\bel{4} (u,\Delta v)_\mu = (\Delta u, v)_\mu \qe for all $u, v \in
\ell^2(V,\mu)$. This follows from the symmetric
weight function, i.e. $w_{ij} = w_{ji}$ for all $i$ and $j$. \\
\noindent Moreover, \item[(ii)] $\Delta$ is nonnegative, i.e.
\bel{5} (\Delta u,u)_\mu \geq 0 \qe for all $u\in \ell^2(V,\mu)$.
This follows from the Cauchy-Schwarz inequality. \item[(iii)]
$\Delta u =0$ iff $u$ is constant.

Clearly, $\Delta u = 0$ if $u$ is constant. Let $\Delta u = 0$ and
assume that $u$ is not constant. Then there exists a vertex, say
$i$, with $u(i) \ge u(j)$ for all $j\sim i$ with strict inequality
for at least one such $j$. Thus there exists a nontrivial local
maximum. This is a contradiction since $\Delta u(i)=0$ implies
that the value $u(i)$ is the average of the values at the
neighbors of $i$. Since $\Gamma$ is connected, $u$ then has to be
a constant. (When $\Gamma$ is not connected, a solution of $\Delta
u =0$ is constant on every connected component of $\Gamma$.)
\end{enumerate}
We say that $\lambda$ is an eigenvalue of $\Delta$ if there exists
some $u\notequiv \ 0$ with \bel{3aa} \Delta u = \lambda u. \qe The
preceding properties have consequences for the eigenvalues of
$\Delta$:
\begin{itemize}
\item By (i), the eigenvalues are real. \item By (ii), they are
nonnegative, i.e. $\lambda_k \geq 0$ for all $k$. \item By (iii),
the smallest eigenvalue is $\lambda_0 =0$. Since $\Gamma$ is
connected, this eigenvalue is simple, i.e. \bel{7} \lambda_k >0
\end{equation}
for $k>0$ where the eigenvalues are ordered as
$$ \lambda_0=0 < \lambda_1 \le ... \le \lambda_{N-1}.$$
\end{itemize}

In the literature, the normalized Laplace operator $\Delta$ is
sometimes introduced in a different way (see for instance
\cite{Sunada}): Let $\mathcal{E}$ be the set of all oriented edges
in $\Gamma$. We denote by $o(e)$ the origin and  by $t(e)$ the
terminus of the edge $e\in\mathcal{E}$, respectively. Furthermore,
$\overline{e}$ denotes the inversion of the edge $e\in
\mathcal{E}$. Let $C^0(V,\mathds{R})= \{u: V \to \mathds{R}\}$ and
$C^1(V,\mathds{R})= \{f: E \to \mathds{R}, f(\bar{e}) = -f(e)\}$
be the 0th and 1st cochain group and let \be \label{Innerprod1}
(u,v)_\mu = \sum_{i\in V}\mu(i) u(i)v(i)\qe and \be
\label{Innerprod2} (f,g)_{\mu}=
\frac{1}{2}\sum_{e\in\mathcal{E}}\mu(e)f(e)g(e)\qe be inner
products on $C^0(V,\mathds{R})$ and $C^1(V,\mathds{R})$,
respectively, where $\mu(i) = d_i$ and $\mu(e) = w_{ij}$ for
$e=(i,j)\in\mathcal{E}$. Note that $C^0(V,\mathds{R})$ together
with the inner product \rf{Innerprod1} is equal to
$\ell^2(V,\mu)$. Now we can define the normalized Laplace operator
as
\[\Delta:C^0(V,\mathds{R})\to C^0(V,\mathds{R}),\]
\be \label{C2} \Delta= d^*d,\qe where $d:C^0(V,\mathds{R}) \to
C^1(V,\mathds{R}),$ \be \label{C88} (dv)(e) = v(t(e)) - v(o(e))\qe
is the coboundary operator ($d$ can be considered as a discrete
analogue of the exterior derivative) and $d^*: C^1(V,\mathds{R})
\to C^0(V,\mathds{R})$
\[(d^*f)(i) = - \frac{1}{\mu(i)}\sum_{e\in \mathcal{E}_i}\mu(e)f(e),\] is the (formal) adjoint of $d$ with respect to the inner
products \rf{Innerprod1} on $C^0(V,\mathds{R})$ and
\rf{Innerprod2} on $C^1(V,\mathds{R})$. Here, $\mathcal{E}_i$  is
given by $\mathcal{E}_i := \{e\in \mathcal{E} : o(e) =i\}$.

A simple calculation shows that \ba \label{C7b} \nonumber
(du,dv)_\mu &=& \frac{1}{2} \sum_{e=(i,j)\in \mathcal{E}}
w_{ij}(u(i)-u(j))(v(i)-v(j))\\&=& \frac{1}{2}\left(\sum_{i,j\in V}
w_{ij}u(i)v(i) + \sum_{i,j\in V}
w_{ij}u(j)v(j) -2 \sum_{i,j\in V}w_{ij}u(i)v(j)\right)\\
\nonumber
&=& \sum_{i\in V} d_i u(i)\left(v(i) - \frac{1}{d_i}\sum_{j\in V}w_{ij}v(j)\right)\\
&=& (u,\Delta v)_\mu.
\end{eqnarray}
This confirms that the two definitions \rf{3} and \rf{C2} of the
normalized Laplace operator coincide.

An orthonormal basis of $\ell^2(V,\mu)$ consisting of
eigenfunctions of $\Delta$,
$$ u_k,\ k=0,...,N-1$$
can be constructed in the standard way which we now recall. Let
$H_0:=H:=\ell^2(V,\mu)$ be the Hilbert space of all real-valued
functions on $\Gamma$ with the inner product $(.,.)_\mu$.  We
iteratively define \bel{7c} H_k:=\{ v\in H: (v,u_i)_\mu=0 \mbox{
for } i\le k-1 \},
\end{equation}
starting with a constant function $u_0$ as the eigenfunction for
the eigenvalue $\lambda_0=0$. Then the $k$th eigenvalue is given
by \bel{7d} \lambda_k= \inf_{u \in
H_k-\{0\}}\frac{(du,du)_\mu}{(u,u)_\mu},
\end{equation}
and the corresponding eigenfunction $u_k$ realizes this infimum.
By way of contrast, the highest eigenvalue is also given by
\bel{LambdaK} \lambda_{N-1}= \sup_{u\;  \smallnotequiv \
0}\frac{(du,du)_\mu}{(u,u)_\mu}.
\end{equation}  In
particular, for any eigenfunction $u$ for some eigenvalue $\lambda
$, we then have \bel{8k} \lambda =\frac{(du,du)_\mu}{(u,u)_\mu}.
\qe All different eigenfunctions are orthogonal to each other. In
particular the eigenfunctions $u_1,\dots u_{N-1}$ are orthogonal
to $u_0$, the eigenfunction for the eigenvalue $\lambda_0=0$. This
implies that \bel{8m} \sum_i d_i u_k(i)=0 \qe for $k=1,\dots
,{N-1}$, since $u_0(i)$ is constant for all $i$.

The largest eigenvalue satisfies \bel{12a} \lambda_{N-1} \le 2 \qe
with equality if and only if $\Gamma$ is bipartite. A
corresponding eigenfunction equals a positive constant $c$ on one
class and $-c$ on the other class of vertices. In contrast, for
loopless graphs, the highest eigenvalue $\lambda_{N-1}$ becomes
smallest on a complete graph $K_N$\footnote{$K_N$ denotes an
unweighted complete loopless graph on $N$ vertices.}, namely
\bel{12b} \lambda_{N-1}= \frac{N}{N-1}.
\end{equation}
By considering the trace of $\Delta$ we obtain \be
\label{trace}\sum_i\lambda_i = N-\sum_i\frac{w_{ii}}{d_i}.\qe
Altogether, the eigenvalues satisfy \bel{9}  0=\lambda_0
<\lambda_1\leq\frac{N-\sum_i\frac{w_{ii}}{d_i}}{N-1}\leq
\lambda_{N-1} \leq 2. \qe Similarly to \rf{12b},  the first
eigenvalue $\lambda_1$ is largest for the complete graph $K_N$,
achieving the bound in \rf{9}, that is \bel{10}
\lambda_1=\frac{N}{N-1}.
\end{equation}
For any other unweighted graph, we have in fact \bel{11}
\lambda_1\le 1. \qe Hence, the complete graph $K_N$ satisfies
\[\lambda_0 = 0 \text{ and }\lambda_1 = \ldots = \lambda_{N-1} = \frac{N}{N-1}.\]
In fact, it is easy to show that for the complete graph $K_N$ any
function $v$ that satisfies $\sum_jd_jv(j)= 0$ is an eigenfunction
for the eigenvalue $\frac{N}{N-1}$.

\section{The Cheeger constant and its dual and eigenvalue
estimates}\label{section3} \label{cheeg} Our starting point are
the estimates for the first eigenvalue $\lambda_1$ in terms of the
(Polya-)Cheeger constant, see for example \cite{Alon85, Alon86,
Dodziuk84, Chung96,Chung97}. The (Polya-) Cheeger constant
\cite{Polya51} of a weighted graph is defined as \bel{ch1}
 h:= \min_U
\frac{|E(U,\overline{U})|}{\min\{\mathrm{vol}(U),
\mathrm{vol}(\overline{U})\}}= \min_{U\subset V :\mbox{\small
vol}(U)\leq  \frac{1}{2}\mbox{\small
vol}(V)}\frac{|E(U,\overline{U})|}{\mathrm{vol}(U)}, \qe where $U$
and $\overline{U}=V\setminus U$ yield a partition of the vertex
set $V$ and $U,\overline{U}$ are both nonempty. Here the volume of
$U$ is given by $\mathrm{vol}(U) := \sum_{i\in U} d_i$,
$E(U,\overline{U})\subseteq E$ is the subset of all edges with one
vertex in $U$ and one vertex in $\overline{U}$, and
$|E(U,\overline{U})|:= \sum_{k\in U, l\in \overline{U}}w_{kl}$ is
the sum of the weights of all edges in $E(U,\overline{U})$. In
general we have $h\leq 1$ and equality holds for instance if
$\Gamma$ is given by $K_2$ or $K_3$. This follows from the
definition of $h$, since $|E(U,\overline{U})| \leq
|E(U,\overline{U})|+ |E(U,U)|= \mathrm{vol}(U)$ and
$|E(U,\overline{U})| \leq |E(U,\overline{U})|+
|E(\overline{U},\overline{U})|= \mathrm{vol}(\overline{U})$.
 Let us first recall \cite{Chung96,Chung97} how
$h$ can bound $\lambda_1$ from above. We use the variational
characterization \rf{7d}, observing that $H_1$ is the set of all
functions $v$ with the normalization $\sum_{i\in V} d_iv(i)=0$.
Let the edge set $E(U,\overline{U})$ divide the graph into the two
disjoint sets $U,\overline{U}$ of nodes, and let $U$ be the one
with the smaller volume $\mathrm{vol}(U)=\sum_{i\in U} d_i$. We
consider a function $v$ that is = 1 on all the nodes in $U$ and
$=-\alpha$ for some positive $\alpha$ on $\overline{U}$. $\alpha$
is chosen so that the normalization $\sum_{i\in V} d_iv(i)=0$
holds, that is, $\sum_{i \in U} d_i - \sum_{i \in \overline{U}}
d_i \alpha =0$. Since $\overline{U}$ is the subset with the larger
volume $\sum_{i\in\overline{U}} d_i$, we have $\alpha \le 1$.
Thus, for our choice of $v$, the quotient in \rf{7d} becomes $\le
\frac{(1+\alpha)^2|E(U,\overline{U})|}{\sum_{i \in U} d_i +\sum_{i
\in \overline{U}} d_i
  \alpha^2}=\frac{(1+\alpha)|E(U,\overline{U})|}{\sum_{i \in U} d_i} \le
2\frac{|E(U,\overline{U})|}{\sum_{i\in U}d_i}=
2\frac{|E(U,\overline{U})|}{\mathrm{vol}(U)}$. Since this holds
for all such splittings of our graph $\Gamma$, we obtain from
\rf{ch1} and \rf{7d} \bel{13a} \lambda_1 \le 2 h.
\end{equation}
As a lower bound for $\lambda_1$ in terms of the Cheeger constant
$h$ we obtain: \bel{13c} \lambda_1 \ge 1-\sqrt{1- h^2}.\qe In
fact, the estimates \rf{13a} and \rf{13c}  hold under rather
general conditions, and an appropriate version is also true for
the algebraic (non-normalized) graph Laplace operator
\cite{Mohar89}.

The crucial step in the proof of \rf{13c} is the next lemma which
we recall here, because we will make use of it in the following.
The proof given here is mainly based on \cite{Diaconis91} and uses
some generalizations that can be found in \cite{Chung97}.
\begin{lemma}\label{55}
Let $g\in \ell^2(V,\mu)$ with $S(g) := \{i \in V : g(i)
>0\}\neq \emptyset$, put
\be \label{hg} h(g) := \min_{\emptyset \neq S \subseteq
S(g)}\frac{|E(S,\overline{S})|} {\mathrm{vol}(S)}\qe
 and let $g_+$ be the positive part of $g$, i.e.
\[g_+(i)=\left\{\begin{array}{ccc} g(i) & \mbox{ if}&  g(i)>0 \\0 &&
    else. \end{array}\right.\]
Then
\[1+\sqrt{1-h^2(g)}\geq \frac{\sum_{e=(i,j)}w_{ij}(g_+(i)-g_+(j))^2}
{\sum_i d_ig_+(i)^2}\geq 1-\sqrt{1-h^2(g)}. \]
\end{lemma}
\begin{proof}For technical reasons it is convenient to define the
new weights $\mu_{ij}$ by $\mu_{ij}=w_{ij}$ for all $i\neq j$ and
$\mu_{ii}=\frac{1}{2}w_{ii}$ for all $i\in V$.

First, we write
\begin{eqnarray*}W&:=& \frac{\sum_{e=(i,j)}w_{ij}(g_+(i)-g_+(j))^2}{\sum_i
d_ig_+(i)^2}\\&=&
\frac{\sum_{e=(i,j)}\mu_{ij}(g_+(i)-g_+(j))^2}{\sum_i d_ig_+(i)^2}\\
&=& \frac{\sum_{e=(i,j)}\mu_{ij}(g_+(i)-g_+(j))^2
\sum_{e=(i,j)}\mu_{ij}(g_+(i)+g_+(j))^2}{\sum_i d_ig_+(i)^2
\sum_{e=(i,j)}\mu_{ij}(g_+(i)+g_+(j))^2} \\&=:&
\frac{I}{II}.\end{eqnarray*} Using the Chauchy-Schwarz inequality
we obtain
\[I \geq \left(\sum_{e=(i,j)}\mu_{ij}|g_+(i)^2-g_+(j)^2|\right)^2
=\left(\sum_{e=(i,j)}w_{ij}|g_+(i)^2-g_+(j)^2|\right)^2.\] Now we
have
\begin{eqnarray*}\sum_{e=(i,j)}w_{ij}|g_+(i)^2-g_+(j)^2|&\geq& \sum_{e=(i,j)}w_{ij}|g_+(i)^2-g_+(j)^2|
\\&=& \sum_{e=(i,j):g_+(i)>
g_+(j)}w_{ij}(g_+(i)^2-g_+(j)^2)\\&=& 2  \sum_{e=(i,j):g_+(i)>
g_+(j)}w_{ij}\int_{g_+(j)}^{g_+(i)}t dt\\&=& 2 \int_{0}^{\infty}
\sum_{e=(i,j):g_+(j)\leq t< g_+(i)}w_{ij} \,tdt.
\end{eqnarray*}Note that $\sum_{e=(i,j):g_+(j)\leq t< g_+(i)}w_{ij}=|E(S_t,\overline{S_t})|$
where $S_t:= \{i:g_+(i)>t\}$. Using \rf{hg} we obtain,
\begin{eqnarray*}\sum_{e=(i,j)}w_{ij}|g_+(i)^2-g_+(j)^2|&\geq& 2
h(g)\int_0^\infty\mathrm{vol}(S_t)tdt\\&=&  2
h(g)\int_0^\infty\sum_{i:g_+(i)>t}d_itdt\\&=&  2 h(g)\sum_{i\in
V}d_i\int_0^{g_+(i)}tdt
\\&=& h(g)\sum_id_ig_+(i)^2
\end{eqnarray*}and so it follows  \[I\geq h^2(g)(\sum_id_ig_+(i)^2)^2.\]
\begin{eqnarray*}II&=& \sum_i d_ig_+(i)^2
\sum_{e=(i,j)}\mu_{ij}(g_+(i)+g_+(j))^2 \\&=& \sum_i d_ig_+(i)^2
(\sum_id_ig_+(i)^2 + \sum_{i,j}w_{ij}g_+(i)g_+(j))\\&=& \sum_i
d_ig_+(i)^2 (2\sum_id_ig_+(i)^2 -
\sum_{e=(i,j)}\mu_{ij}(g_+(i)-g_+(j))^2)\\&=& (2-W)(\sum_i
d_ig_+(i)^2)^2.
\end{eqnarray*}
Combining everything we obtain,
\[W\geq \frac{h^2(g)}{(2-W)}\] and  consequently
\[1+\sqrt{1-h^2(g)}\geq W \geq 1-\sqrt{1-h^2(g)}.\]
\end{proof}
The second observation that we need to prove the Cheeger
inequality \rf{13c} is the following lemma \cite{Diaconis91}:
\begin{lemma}
For every non-negative real number $\lambda$ and
$g\in\ell^2(V,\mu)$ we have
\[\lambda \geq \frac{\sum_{e=(i,j)}w_{ij}(g_+(i)-g_+(j))^2}{\sum_i
d_ig_+(i)^2}=W \] if $\Delta g(i)\leq \lambda g(i)$ for all  $i\in
S(g)$.
\end{lemma}
\begin{proof}We have
\[(\Delta g, g_+)_\mu = \sum_{i\in V}d_i\Delta g(i)g_+(i)
\leq \lambda\sum_{i\in S(g)}d_i g_+(i)g_+(i) = \lambda\sum_{i\in
V}d_i g_+(i)g_+(i)\] and \begin{eqnarray*}(\Delta g, g_+)_\mu &=&
(d g, d g_+)_\mu = \sum_{e=(i,j)\in
E}w_{ij}(g(i)-g(j))(g_+(i)-g_+(j))
\\&\geq&\sum_{e=(i,j)\in E}w_{ij}(g_+(i)-g_+(j))^2.
\end{eqnarray*}
\end{proof}
The Cheeger inequality now follows from the last two lemmata by
taking $\lambda=\lambda_1$ and $g=u_1$ an eigenfunction for
$\lambda_1$. Since $(u_1,\mathds{1})_\mu=0$ we have $S(u_1)\neq
\emptyset$ and $T(u_1):= \{i\in V: u_1(i) <0\}\neq \emptyset$.
This implies that it is always possible to choose $u_1$ such that
$\mathrm{vol}(S(u_1))\leq \mathrm{vol}(\overline{S(u_1)})$ (if
$\mathrm{vol}(S(u_1))\geq \mathrm{vol}(\overline{S(u_1)}$) take
$-u_1$ instead of $u_1$) and thus $h(u_1)\geq h$.

In any case, in qualitative terms, the Cheeger inequalities
\rf{13a} and \rf{13c}  simply say that $\lambda_1$ becomes small
when the graph can be easily (that is, by cutting only few edges)
decomposed into two large parts. Thus, $\lambda_1$ is small, that
is, close to its minimal value $0$, when $\Gamma$ is similar to a
disconnected graph, with equality iff $\Gamma$ is disconnected
itself. Similarly, and this brings us to our topic, the largest
eigenvalue is large, that is, close to its maximal value $2$, when
$\Gamma$ is close to a bipartite graph, with equality iff
$\Gamma$ is bipartite itself.  \\
The main purpose of this section then is a dual version of
\rf{13a} and \rf{13c} for the largest eigenvalue $\lambda_{N-1}$.
More precisely, we shall obtain an estimate for $\lambda_{N-1}$ in
terms of a dual version of the Cheeger constant which we now
introduce. Let $V_1, V_2$ and $\overline{V_1\cup V_2} =: V_3$ be a
partition of the vertex set $V$ into three disjoint sets such that
$V_1$ and $V_2$ are nonempty.

It is helpful to think of $V_1\cup V_2$ as the (almost) bipartite
part of $\Gamma$  and $V_3$ as the part of $\Gamma$ that contains
many cycles of odd length, i.e. $V_3$ is not bipartite.

For a partition $V_1,V_2,V_3$ of the vertex set $V$ we define:
\be\label{hh}\overline{h}:= \max_{V_1,V_2}
\frac{2|E(V_1,V_2)|}{\mathrm{vol}(V_1)+\mathrm{vol}(V_2)},\qe
where as before the volume of $V_k$ is given by $\mathrm{vol}(V_k)
:= \sum_{i\in V_k} d_i$ and $|E(V_i,V_j)|:= \sum_{k\in V_i, l\in
V_j}w_{kl}$.

The next theorem shows that $\overline{h}$ characterizes bipartite
graphs.
\begin{theo}\label{Lemma2} $\overline{h} \le 1$, and $\overline{h} = 1$ if and only if $\Gamma$ is
bipartite.
\end{theo}
\begin{proof}First, note that, for a partition $V_1,V_2$ and $V_3$ of $V$, the volume
of $V_i$ can also be written in the form
\begin{equation}\label{71}
\mathrm{vol}(V_i) = \sum_{j=1}^3|E(V_i,V_j)|
\end{equation}
Consequently, $\overline{h}$ is given by
\begin{equation}\label{A4} \overline{h} = \max_{V_1,V_2}
\frac{2|E(V_1,V_2)|}{\sum_{j=1}^3|E(V_1,V_j)| +
\sum_{j=1}^3|E(V_2,V_j)| }.\end{equation} Thus, clearly, \be
\label{5h} \overline{h}\leq1.\qe

Assume that $\Gamma$ is bipartite. Then there exists a partition
$V_1,V_2,V_3$ of $V$ such that $V_3=\emptyset$ and there are no
edges within the subsets $V_1$ and $V_2$. For this partition it
follows that $|E(V_1,V_3)| = |E(V_2,V_3)| = |E(V_1,V_1)| =
|E(V_2,V_2)| = 0$. By \rf{A4} this implies that $\overline{h}
\geq1$. Together with \rf{5h} it follows that $\overline{h} =1$.

Now assume that $\overline{h} = 1$. Equation \rf{A4} implies that
there exists a partition $V_1,V_2,V_3$ of $V$ such that
$|E(V_1,V_3)| = |E(V_2,V_3)| = |E(V_1,V_1)| = |E(V_2,V_2)| = 0$.
Since $\Gamma$ is connected $V_3 = \emptyset$ and thus $\Gamma$ is
bipartite.
\end{proof}
As an illustration, let us consider loopless Erd\"os-Renyi random
graphs, i.e. we start with a given vertex set and add edges
between two vertices with a fixed probability $p$. If we start
with $p=1$ then we obtain a complete graph and thus
$\bar{h}\approx1/2$, as will be shown in Example \ref{23}. Now if
we decrease $p$ we decrease the number of edges in the graph. This
will lead to a local bipartite subgraph in $\Gamma$ and thus
$\bar{h}$ will be increased. If we decrease $p$ further we finally
have $|E|\approx|V|$, i.e. the graph will be approximately a tree
(we assume that the random graphs are connected) thus $\bar{h}
\approx 1$. We conclude that, for random graphs $\bar{h}$ is a
function of $p$. More details are revealed by numerical
simulations.
\begin{satz}\label{satz}
For a loopless graph $\Gamma$,  \[\frac{1}{2}\leq \overline{h}.\]
\end{satz}
\begin{proof}Assume that there exists a
partition $V_1,V_2$ and $V_3 = \emptyset$ of the vertex set $V$
such that \be \label{24} |E(V_1,V_2)| \geq
\max_{i=1,2}|E(V_i,V_i)|.\qe Then by using \rf{71} we obtain:
\begin{eqnarray*}\overline{h}&\geq& \max_{V_1,V_2,
V_3=\emptyset}\frac{2|E(V_1,V_2)| }{2|E(V_1,V_2)|
+|E(V_1,V_1)|+|E(V_2,V_2)| }\\ &\geq& \max_{V_1,V_2,
V_3=\emptyset}\frac{|E(V_1,V_2)| }{|E(V_1,V_2)|
+\max_{i=1,2}|E(V_i,V_i)|}\\&\geq&\frac{1}{2}.\end{eqnarray*}
Thus, it is sufficient to find a partition that satisfies \rf{24}.

In the following we will construct such a partition. Start with an
arbitrarily partition $V_1,V_2$ and $V_3 = \emptyset$ of $V$. If
\rf{24} is satisfied we are done. Otherwise, assume w.l.o.g. that
$|E(V_1,V_1)|
> |E(V_1,V_2)|$, i.e. \[\sum_{i\in V_1}\sum_{j\in V_1}w_{ij}>
\sum_{i\in V_1}\sum_{j\in V_2}w_{ij}.\] We observe that there
exists a vertex $i$ in $V_1$ such that
\[\sum_{j\in V_1}w_{ij} > \sum_{j\in V_2}w_{ij}.\] We remove the
vertex $i$ from $V_1$ and add it to $V_2$. By doing so,
$|E(V_1,V_2)|$ is increased by $\sum_{j\in V_1}w_{ij}- \sum_{j\in
V_2}w_{ij}>0$, $|E(V_1,V_1)|$ is decreased by $\sum_{j\in
V_1}w_{ij}$, and $|E(V_2,V_2)|$ is increased by $\sum_{j\in
V_2}w_{ij}$.  If \rf{24} is still not satisfied, continue this
procedure several times. Eventually, \rf{24} holds since
$|E(V_1,V_2)|$ is strictly monotonically increasing.
\end{proof}
This lower bound is optimal, since Example \ref{23} shows that for
complete graphs $K_N$, $\overline{h}(K_N)\rightarrow\frac{1}{2}$
as $N\rightarrow \infty$. Clearly, the proof of Proposition
\ref{satz} cannot be extended to graphs with loops. In fact,
Proposition \ref{satz} only holds for loopless graphs, as can be
seen by considering a graph with $\sum_i\frac{w_{ii}}{d_i}>N-1$.
In that case, $\overline{h}\geq 1/2$ would lead to a contradiction
in Theorem \ref{A6} since
\[1\leq 2\overline{h} \leq \lambda_{N-1} < 1, \] where we used \rf{trace} in the last inequality.

We now have a counterpart of the Cheeger inequality \rf{13a} and
\rf{13c}.
\begin{theo}\label{A6}
The largest eigenvalue $\lambda_{N-1}$ of the graph Laplace
operator $\Delta$ satisfies  \begin{equation}\label{A2}
2\overline{h} \leq \lambda_{N-1} \leq
1+\sqrt{1-(1-\overline{h})^2}.\end{equation}
\end{theo}
\begin{proof}First, we prove that $2\overline{h} \leq \lambda_{N-1}$.
The largest eigenvalue $\lambda_{N-1}$ of $\Delta$ is given by
\rf{LambdaK}. Let $V_1$, $V_2$, $V_3$ be a partition that achieves
$\overline{h}$. We consider the following function $u$:
\[u(i) = \left\{\begin{array}{lll}\frac{1}{\mathrm{vol}(V_1)} & \mbox{if} & i \in V_1 \\
\frac{-1}{\mathrm{vol}(V_2)} & \mbox{if} & i \in V_2 \\ 0 &&
\mbox{else.}
\end{array} \right.\]
Substituting $u$ in (\ref{LambdaK}) yields:
\begin{eqnarray*}
\lambda_{N-1}&=& \sup_{v\neq 0}\frac{(dv,dv)_\mu}{(v,v)_\mu} =
\sup_{v\neq0}\frac{\sum_{e=(i,j)}w_{ij}(v(i)-v(j))^2}{\sum_id_iv(i)^2}\\
&\geq& \frac{\left( \frac{1}{\mathrm{vol}(V_1)} +
\frac{1}{\mathrm{vol}(V_2)}\right)^2|E(V_1,V_2)| +
\left(\frac{1}{\mathrm{vol}(V_1)}\right)^2 |E(V_1, V_3)|+
\left(\frac{1}{\mathrm{vol}(V_2)}\right)^2 |E(V_2, V_3)| }{\left(
\frac{1}{\mathrm{vol}(V_1)} + \frac{1}{\mathrm{vol}(V_2)}\right)}
\\&\geq& \frac{(\mathrm{vol}(V_1)+\mathrm{vol}(V_2))^2}{2\mathrm{vol}(V_1)\mathrm{vol}(V_2)}
\frac{2|E(V_1,V_2)|}{\mathrm{vol}(V_1)+ \mathrm{vol}(V_2)} +
\frac{\min(\mathrm{vol}(V_1),\mathrm{vol}(V_2))}{\max(\mathrm{vol}(V_1),\mathrm{vol}(V_2))}\frac{|E(V_1\cup
V_2, V_3)|}{(\mathrm{vol}(V_1) + \mathrm{vol}(V_2))}
 \\&\geq& 2 \overline{h}+
\frac{\min(\mathrm{vol}(V_1),\mathrm{vol}(V_2))}{\max(\mathrm{vol}(V_1),\mathrm{vol}(V_2))}\frac{|E(V_1\cup
V_2, V_3)|}{(\mathrm{vol}(V_1) + \mathrm{vol}(V_2))}
\\&\geq& 2\overline{h},
\end{eqnarray*}where we used the simple inequality
$\frac{(a+b)^2}{2ab}\geq 2$ for $a,b\in\R$.

Now we prove the remaining inequality $\lambda_{N-1} \leq
1+\sqrt{1-(1-\overline{h})^2}$. When one studies the largest
eigenvalue of $\Delta$ it is convenient to introduce the operator
$L=2I -\Delta$. If $\lambda$ is an eigenvalue of $\Delta$ and
corresponding eigenfunction $u$ then $u$ is also an eigenfunction
for $L$ and corresponding eigenvalue $\mu = 2-\lambda$. Thus,
controlling the largest eigenvalue $\lambda_{N-1}$ of $\Delta$
from above is equivalent to controlling the smallest eigenvalue
$\mu_0$ of $L$ from below. The smallest eigenvalue $\mu_0$ of $L$
is given by
\begin{eqnarray*}\mu_0 &=& \inf_{u\neq 0} \frac{\frac{1}{2}\sum_{i,j\in V}w_{ij}(u(i) +
u(j))^2}{\sum_{i\in V}d_iu(i)^2}\\&=&\inf_{u\neq 0}
\frac{\sum_{e=(i,j)\in E}\mu_{ij}(u(i) + u(j))^2}{\sum_{i\in
V}d_iu(i)^2}
\end{eqnarray*} where as above $\mu_{ij}=w_{ij}$ for all $i\neq j$ and
$\mu_{ii}=\frac{1}{2}w_{ii}$ for all $i\in V$. This simply follows
from the standard minmax characterization of eigenvalues
\[\mu_0=\inf_{u\neq 0}\frac{(Lu,u)_\mu}{(u,u)_\mu}.\] We have for all
$u,v\in\ell^2(V,\mu)$
\begin{eqnarray*}
(Lu,v)_\mu&=& \sum_id_iLu(i)v(i)= \sum_i\sum_jw_{ij}
(u(i)+u(j))v(i)\\&=& \sum_j\sum_iw_{ji} (u(j)+u(i))v(j)
\end{eqnarray*} where we just exchanged $i$ and $j$. Adding the last two
lines and setting $u=v$ yields\[(Lu,u)_\mu=
\frac{1}{2}\sum_{i,j}w_{ij}(u(i)+u(j))^2.\]

In order to prove the lower bound for $\mu_0$ we will use a
technique developed in \cite{Desai94}. The idea is the following:
Construct a graph $\Gamma^\prime$ out of $\Gamma$ s.t. the
quantity $h^\prime(g)$ defined  in Lemma \ref{55} for the new
graph $\Gamma^\prime$ controls $\mu_0$ from below. In a second
step, we show that $h^\prime(g)$ in turn can be controlled by the
quantity $1-\overline{h}$ of the original graph. This then yields
the desired estimate.

Let $u$ be an eigenfunction for the eigenvalue $\mu_0$ and define
as above $S(u)= \{i\in V : u(i)>0\}$ and $T(u)= \{i\in V : u(i)
<0\}$. Since $u$ is also an eigenfunction for $\lambda_{N-1}$ of
$\Delta$ we know that $(u,\mathds{1})_\mu=0$ and thus
$S(u),T(u)\neq \emptyset$. Then the new graph $\Gamma^\prime
=(V^\prime, E^\prime)$ is constructed from $\Gamma$ in the
following way. Duplicate all vertices in $S(u)\cup T(u)$ and
denote the copies by a prime, e.g. if $i\in S(u)$ then the copy of
$i$ is denoted by $i^\prime$. The copies of $S(u)$ and $T(u)$ are
denoted by $S^\prime(u)$ and $T^\prime(u)$ respectively. The
vertex set $V^\prime $ of $\Gamma^\prime$ is given by $V^\prime=
V\cup S^\prime(u)\cup T^\prime(u)$. Every edge $(i,j)\in
E(S(u),S(u))$ in $\Gamma$ is replaced by two edges $(i,j^\prime)$
and $(j,i^\prime)$ in $\Gamma^\prime$ s.t. $w_{ij} =
w^\prime_{ij^\prime}=w^\prime_{ji^\prime}$. Similarly, if the edge
is a loop, then $e=(i,i)$ is replaced by one edge $(i,i^\prime)$
s.t. $w_{ii}=w_{ii^\prime}$. The same is done with edges in
$E(T(u),T(u))$. All other edges remain unchanged, i.e. if
$(k,l)\in E\setminus (E(S(u),S(u))\cup E(T(u),T(u)))$ then
$(k,l)\in E^\prime$ and $w_{kl} = w^\prime_{kl}$. It is important
to note that this construction does not change the degrees of the
vertices in $V^\prime\setminus (S^\prime(u)\cup T^\prime(u))$.

Consider the function $g:V^\prime \rightarrow \R$,
\[g(i)=\left\{\begin{array}{ccc} |u(i)| & \mbox{ if}&  i \in S(u)\cup T(u) \\0 &&
else. \end{array}\right.\] It can easily be checked that by
construction of $\Gamma^\prime$ we have
\begin{eqnarray*}\mu_0 &=& \frac{\sum_{e=(i,j)\in E}\mu_{ij}(u(i) +
u(j))^2}{\sum_{i\in V}d_iu(i)^2}\\&\geq&
\frac{\sum_{e^\prime=(i,j)\in E^\prime}w^\prime_{ij}(g(i) -
g(j))^2}{\sum_{i\in V^\prime}d^\prime_ig(i)^2}\\&\geq& 1-
\sqrt{1-(h^\prime(g))^2}
\end{eqnarray*} where we used  Lemma \ref{55} to obtain the last
inequality. For any non-empty subset $W\subseteq S(g) = S(u)\cup
T(u)$ we define $S_1=W\cap S(u)$ and $T_1=W\cap T(u)$. Let
$\emptyset\neq U\subseteq S(g)$ the subset that realizes the
infimum, i.e.
\begin{eqnarray*}h^\prime(g) &=& \inf_{\emptyset\neq W\subseteq S(g)}
\frac{|E^\prime(W,\overline{W})|}{\mathrm{vol}(W)}=
\frac{|E^\prime(U,\overline{U})|}{\mathrm{vol}(U)}\\&=&
\frac{|E(S_1,S_1)|+ |E(T_1,T_1)| + |E(S_1\cup
T_1,\overline{S_1\cup T_1})|}
{\mathrm{vol}(S_1)+{\mathrm{vol}(T_1)}}\\&=& 1- \frac{2
|E(S_1,T_1)|}{\mathrm{vol}(S_1)+{\mathrm{vol}(T_1)}}\\&\geq&1-\overline{h}.
\end{eqnarray*}
 Thus we have
\[2-\lambda_{N-1}=\mu_0 \geq 1 - \sqrt{1-(1-\overline{h})^2}\] and so
\[\lambda_{N-1}\leq 1+\sqrt{1-(1-\overline{h})^2}.\]
\end{proof}

For example, the lower estimate for $\lambda_{N-1}$ in \rf{A2} is
sharp if $\Gamma$ is a bipartite (by Theorem \ref{Lemma2}) or if
$\Gamma$ is a complete graph $K_N$, and $N$ is even (by Example
\ref{23}). In both examples the partition that achieves
$\overline{h}$ satisfies $V_3 =\emptyset$. In fact, the proof of
Theorem \ref{A6} shows that the estimate for $\lambda_{N-1}$ from
below  can only be sharp if $V_3 =\emptyset$. However, if the
volume of $V_3$ is sufficiently large, we can improve the estimate
given in \rf{A2} and estimate the eigenvalue $\lambda_{N-1}$ from
below by using both the Cheeger constant $h$ and its dual
$\overline{h}$.
\begin{coro}
Assume that $V_1,V_2$ and $V_3$ is a partition of $V$ that
achieves $\overline{h}$. If \[\mathrm{vol}(V_1\cup V_2) \leq
\mathrm{vol}(V_3),\] then \be \lambda_{N-1} \geq 2\overline{h} +
\mathcal{R}(V_1,V_2) h,\qe where we define for the partition
$V_1,V_2, V_3$ of the vertex set $V$ \be \label{31}
\mathcal{R}(V_1,V_2) :=
\frac{\min(\mathrm{vol}(V_1),\mathrm{vol}(V_2))}{\max(\mathrm{vol}(V_1),\mathrm{vol}(V_2))}.
\qe
\end{coro}
\begin{proof}The proof of Theorem \ref{A6} shows that
\begin{eqnarray*}
\lambda_{N-1} &\geq& 2 \overline{h}+
\mathcal{R}(V_1,V_2)\frac{|E(V_1\cup V_2, V_3)|}{\mathrm{vol}(V_1)
+ \mathrm{vol}(V_2)} \\&=& 2 \overline{h}+
\mathcal{R}(V_1,V_2)\frac{|E(V_1\cup V_2,
V_3)|}{\min(\mathrm{vol}(V_1\cup V_2) ,
\mathrm{vol}(V_3))}\\&\geq& 2\overline{h} + \mathcal{R}(V_1,V_2)
h.
\end{eqnarray*}
\end{proof}

The next corollary shows that if the eigenfunction for the largest
eigenvalue $\lambda_{N-1}$ is sufficiently localized, then
$\lambda_{N-1}$ can also be controlled from above in terms of the
Cheeger constant $h$
\begin{coro}\label{62}
Let $u$ be the eigenfunction for the largest eigenvalue of
$\Delta$. If the eigenfunction is sufficiently localized, i.e. \be
\label{72} \sum_{i:u(i)\neq0}d_i \leq \sum_{i:u(i)=0}d_i \qe then
\[\lambda_{N-1}\leq 1+\sqrt{1-h^2}.\]
\end{coro}
\begin{proof}
Again we consider the smallest eigenvalue $\mu_0$ of the operator
$L=I+P$ instead of the largest eigenvalue $\lambda_{N-1}$ of
$\Delta$. Since $u$ is also an eigenfunction for the eigenvalue
$\mu_0$ we have
\begin{eqnarray*}
\mu_0 &=&  \frac{\sum_{e=(i,j)}\mu_{ij}(u(i) +
u(j))^2}{\sum_{i}d_iu(i)^2}\\&\geq&
\frac{\sum_{e=(i,j)}w_{ij}||u(i)|-
|u(j)||^2}{\sum_{i}d_i|u(i)|^2}\\&\geq& 1-\sqrt{1-h^2(|u|)}
\end{eqnarray*}where we used the reverse triangle inequality and Lemma
\ref{55}. Since the eigenfunction is sufficiently localized, it
follows from \rf{72} that  $\mathrm{vol}(S(u)\cup T(u)) \leq
\mathrm{vol}(\overline{S(u)\cup T(u)})$. This implies that
$h(|u|)$ satisfies $h(|u|)\geq h$. Inserting this in the above
equation yields \[2-\lambda_{n-1} = \mu_0 \geq 1-\sqrt{1-h^2}.\]
\end{proof}By considering a bipartite graph it becomes clear that
Corollary \ref{62} is in general not true if \rf{72} is not
satisfied.  However, for infinite graphs the situation is
different. For infinite graphs, the supremum of the spectrum can
always be controlled from above in terms of a Cheeger constant,
see Theorem \ref{Fuji}. In Section \ref{section7} we develop a new
approach that allows us to control the largest eigenvalue of a
graph in terms of the Cheeger constant of its neighborhood graphs
(instead of the Cheeger constant of the graph itself), see
Corollary \ref{coro5}.
\begin{satz}\label{56}Let $u$ be an eigenfunction for the largest eigenvalue of
$\Delta$ that satisfies $\max_i|u(i)| =1$ then
\[\lambda_{N-1}\leq 2- \frac{\min_{i,j}w_{ij}\left( 1 -\min_i|u(i)|
\right)^2}{D\mbox{ vol}(V)}\] where $D$ is the diameter of the
graph. In particular, if there exists a vertex $i$ such that $u(i)
=0$ then
\[\lambda_{N-1}\leq 2- \frac{\min_{i,j}w_{ij}}{D\mbox{
vol}(V)}.\]
\end{satz}
\begin{proof}Again, we consider the smallest eigenvalue $\mu_0$ of
$L=I+P$ instead of the largest eigenvalue $\lambda_{N-1}$ of
$\Delta$. Let $i_k, k=1,\ldots,n$ be the shortest path connecting
the vertices that satisfy $\max_i|u(i)| = |u(i_1)|=1$ and
$\min_i|u(i)| = |u(i_n)|$.  Then we have
\begin{eqnarray*} \mu_0 &=&  \frac{\sum_{e=(i,j)}\mu_{ij}(u(i) +
u(j))^2}{\sum_{i}d_iu(i)^2}\\&\geq&
\frac{\sum_{e=(i,j)}w_{ij}||u(i)|-
|u(j)||^2}{\sum_{i}d_i|u(i)|^2}\\&\geq& \min_{i,j}w_{ij}
\frac{\sum_{k=1}^n(|u(i_k)|- |u(i_{k+1})|)^2}{\mbox{ vol}(V)}.
\end{eqnarray*}
Using the Cauchy-Schwarz inequality we obtain
\begin{eqnarray*}
\mu_0&\geq&
\frac{\min_{i,j}w_{ij}}{n}\frac{\left(\sum_{k=1}^n|u(i_k)|-
|u(i_{k+1})|\right)^2}{\mbox{ vol}(V)}\\&\geq&
\frac{\min_{i,j}w_{ij}}{D}\frac{\left( 1 -\min_i|u(i)|
\right)^2}{\mbox{ vol}(V)}
\end{eqnarray*}where we used the
fact that the length of a shortest path connecting any two
vertices is less or equal to $D$. Since $2-\lambda_{N-1} = \mu_0$
the proof is complete.
\end{proof}
In particular, the estimate in Proposition \ref{56} is sharp for
bipartite graphs, because then $|u(i)| = 1$ for all $i \in V$.

Again, by using the concept developed in Section \ref{section7} we
can derive a similar result to Proposition \ref{56}. In Corollary
\ref{73}, we show that the largest eigenvalue can, independently
of the corresponding eigenfunction, be controlled from above in
terms of the diameter of the neighborhood graph.

Jerrum and Sinclair have shown how one can bound the Cheeger
constant $h$ by using canonical paths \cite{Sinclair89,
Sinclair92}. Similarly, we can derive an upper bound for the dual
Cheeger constant $\overline{h}$ by considering a suitable
collection of paths. Let $\sigma_i$ be a path from vertex $i$ to
vertex $i$ with an odd number of edges and let $\Sigma$ be the
collection of all these paths (one for each vertex).
\begin{theo}\label{TheoXi}
We have \[\overline{h} \leq 1-\frac{1}{\xi},\] where \[\xi :=
\max_{e=(k,l)} \frac{1}{w_{kl}} \sum_{i:\sigma_i\ni
e=(k,l)}d_i.\]The sum is over all $i$ for which the path
$\sigma_i$ contains the edge $e=(k,l)$.
\end{theo}
\begin{proof}For simplicity we define the subset $\Omega \subset E$
as $\Omega := E(V_1,V_1)\cup E(V_1,V_3)\cup E(V_2,V_2)\cup
E(V_2,V_3)$. Now observe that for every vertex $i \in V_1\cup V_2$
a path $\sigma_i$ with an odd number of edges contains at least
one edge in $\Omega$. Thus we have for any partition $V_1,V_2,V_3$
of the vertex set $V$
\begin{eqnarray*} \mathrm{vol}(V_1) +\mathrm{vol}(V_2)&=& \sum_{i\in
V_1\cup V_2}d_i\\&\leq& \sum_{e=(k,l)\in \Omega}
\sum_{i:\sigma_i\ni e=(k,l), i\in V_1\cup V_2}d_i\\&\leq&
\sum_{e=(k,l)\in \Omega} \sum_{i:\sigma_i\ni
e=(k,l)}\frac{w_{kl}}{w_{kl}}d_i\\&\leq& \xi
\sum_{e=(k,l)\in\Omega}w_{kl} \leq \xi|\Omega|.
\end{eqnarray*}
 Since this holds for all partitions, we have for the partition $V_1,V_2$ and $V_3$ that achieves $\overline{h}$
\[1-\overline{h} =
\frac{|\Omega|}{\mathrm{vol}(V_1) +\mathrm{vol}(V_2)}\geq
\frac{1}{\xi}.\]
\end{proof}
\begin{coro}We have
\[\xi \leq d_\Gamma w_\Gamma b_\Gamma,\]
where $d_\Gamma = \max_i d_i,$ $w_\Gamma =
\frac{1}{\min_{i,j}w_{ij}}$ and $b_\Gamma = \max_e\#\{\sigma\in
\Sigma: e\in \sigma\}$. Together with Theorem \ref{A6} and Theorem
\ref{TheoXi}  this implies that \be\label{57} \lambda_{N-1}\leq
1+\sqrt{1-\left(\frac{1}{d_\Gamma w_\Gamma b_\Gamma}\right)^2}.\qe
\end{coro}
\begin{rem}
Diaconis and Stroock show in \cite{Diaconis91}, by using a
discrete analog of the Poincare inequality, that the largest
eigenvalue satisfies \be\label{58}\lambda_{N-1}\leq
2-\frac{2}{d_\Gamma w_\Gamma b_\Gamma \sigma_\Gamma}\qe where
$\sigma_\Gamma$ is the maximum number of edges in any
$\sigma\in\Sigma$. A simple calculation shows that the estimate
\rf{57} obtained from the dual Cheeger inequality is better than
the estimate \rf{58} obtained from the Poincare inequality iff
\[d_\Gamma w_\Gamma b_\Gamma < \frac{1}{\sigma_\Gamma} + \frac{\sigma_\Gamma}{4}.\]
In general, it is not clear which of these estimates is better.
See also the related discussion in \cite{Fulman99}, where the
authors analyze when the Cheeger inequality improves the Poincare
estimate for the smallest nontrivial eigenvalue.
\end{rem}

\section{Relations between $h$ and $\overline{h}$ }

By looking at the definitions of $h$ and $\overline{h}$ it is
apparent that there is a connection between those two quantities.
We shall explore this now in more detail.

Similarly to \rf{31} we define:
\begin{defi}For any partition $U,\overline{U}$ of the
vertex set $V$ we define
\[\mathcal{R}(U) :=\frac{\min(\mathrm{vol}(U),
\mathrm{vol}(\overline{U}))}{\max(\mathrm{vol}(U),
\mathrm{vol}(\overline{U}))}.\] Furthermore,
\[ \mathcal{R} :=\max_U \mathcal{R}(U).\]
\end{defi}
First, we will restrict ourselves to unweighted graphs. Later on
we will prove similar results for weighted graphs. \\\\
\textbf{Unweighted graphs:}\\
\begin{lemma}\label{16} Let $\Gamma$ be an unweighted graph with $N$ vertices, then
\begin{equation}
\frac{N-1}{N+1}\leq \mathcal{R}\leq 1.
\end{equation}Equality holds on the left hand side if and only if $\Gamma$ is a regular
graph and $N$ is odd.
\end{lemma}
\begin{proof}
Order the vertices w.r.t. their degree, i.e. $d_1 \geq d_2\geq
\ldots \geq d_N$.  We construct a partition $U,\overline{U}$ of
$V$ that satisfies $\frac{N-1}{N+1} \leq \mathcal{R}(U)$. We begin
with two empty sets $U_0, \overline{U}_0$. After the partition of
$K$ vertices we denote the subsets by $U_K, \overline{U}_K$.
Having started with vertex $1$ as one of  largest degree, we
iteratively partition the vertices into two subsets such that
vertex $K+1$ is then added to the subset $U_K$, $\overline{U}_K $
that has the smaller volume. We continue this procedure until we
obtain a complete partition $U, \overline{U} := U_N,
\overline{U}_N$ of the vertex set $V$. Let $M\leq N$ be such that
\be\label{defM} \mathrm{vol}(U_{M-1})
\geq\mathrm{vol}(\overline{U}_{M-1})\qe and
\[\mathrm{vol}(\overline{U}_{K}) \geq\mathrm{vol}(U_{K})
\quad \mbox{for }  M\leq K\leq N.\] For simplicity we
define\[\mathrm{vol}(U_K) + \mathrm{vol}(\overline{U}_K)=:
\mathrm{vol}(V_K)\, \mbox{for }1\leq K\leq N.\] Then we have \be
\label{28} \mathrm{vol}(\overline{U}_M)- \mathrm{vol}(U_M) \leq
d_M\leq \frac{\mathrm{vol}(U_M) + \mathrm{vol}(\overline{U}_M)}{M}
= \frac{\mathrm{vol}(V_M)}{M}.\qe Equality holds on the left hand
side if and only if vol$(U_{M-1}) =$ vol$(\overline{U}_{M-1})$ and
equality holds on the right hand side if and only if $d_1 = d_2 =
\ldots = d_M$. For the final partition $U_N,\overline{U}_N$ we
obtain
\begin{eqnarray*}\mathrm{vol}(\overline{U}_N)- \mathrm{vol}(U_N) &=&
\mathrm{vol}(\overline{U}_M) -\left(\mathrm{vol}(U_M) +
\mathrm{vol}(V_N)- \mathrm{vol}(V_M) \right)
\\
&\leq& \frac{\mathrm{vol}(V_M)}{M}-\mathrm{vol}(V_N) +
\mathrm{vol}(V_M)\\&\leq&\frac{\mathrm{vol}(V_N)}{N}.\end{eqnarray*}
The last inequality follows from
\begin{eqnarray}&&\nonumber
\frac{\mathrm{vol}(V_N)}{N} +\mathrm{vol}(V_N) -
\mathrm{vol}(V_M)- \frac{\mathrm{vol}(V_M)}{M}\\&=& \nonumber
\frac{1}{NM}\left[ M(1+N)\mathrm{vol}(V_N) -
N(1+M)\mathrm{vol}(V_M)\right]\\&=&\nonumber
\frac{1}{NM}\left[M(1+N)\left(\mathrm{vol}(V_N)-\mathrm{vol}(V_M)
\right) + (M-N) \mathrm{vol}(V_M)\right]\\&\geq& \label{29}
\frac{1}{NM}\left[M(1+N)(N-M) +(M-N)MN \right]\\&=&\nonumber
\frac{(N-M)}{N} \geq0
\end{eqnarray}
Thus, we constructed a partition that satisfies
\[\max(\mathrm{vol}(U_N), \mathrm{vol}(\overline{U}_N))-
\min(\mathrm{vol}(U_N), \mathrm{vol}(\overline{U}_N))\leq
\frac{\mathrm{vol}(V_N)}{N}.\]Since
\[\frac{\mathrm{vol}(V_N)}{N}=\frac{\max(\mathrm{vol}(U_N),
\mathrm{vol}(\overline{U}_N))+ \min(\mathrm{vol}(U_N),
\mathrm{vol}(\overline{U}_N))}{N},\] this yields
\[\frac{N-1}{N+1} \leq \mathcal{R}(U)\leq \mathcal{R}.\]
From the proof we see that equality holds iff $N=M$, the graph is
regular, and vol$(U_{M-1}) =$ vol$(\overline{U}_{M-1})$. This
implies that equality holds iff $\Gamma$ is regular and $N$ is
odd.
\end{proof}

The last lemma shows that for large graphs, i.e. $N$ large, it is
always possible to partition $V$ into two subsets of almost equal
volume.
\begin{coro}\label{19} In particular we have for unweighted graphs,
\[\frac{1}{2}\leq \mathcal{R}\leq 1\] and equality holds on the left hand side iff $\Gamma$ is
a triangle.
\end{coro}
\begin{proof}The proof follows from Lemma \ref{16} since there
exists only one connected graph on $2$ vertices for which we have
$\mathcal{R} =1$ and the only regular graph on $3$ vertices is the
triangle.
\end{proof}
\begin{theo}\label{17}For unweighted graphs we have
\bel{11}\frac{N-1}{N} h\leq \frac{2\mathcal{R}}{1+\mathcal{R}}
h\leq\overline{h}\leq 1 .\qe Equality holds on the l.h.s. iff
$\Gamma$ is a regular graph and $N$ is odd.
\end{theo}
\begin{proof}
 By \rf{hh}, we have for any partition $U, \overline{U}$ of the vertex set $V$
 \be h\leq
\frac{|E(U,\overline{U})|}{\min(\mathrm{vol}(U),\mathrm{vol}(\overline{U}))}.
\qe For the partition $V_1 = U$, $V_2 = \overline{U}$ and $V_3 =
\emptyset$ we obtain \be \overline{h} \geq
\frac{2|E(U,\overline{U})|}{\mathrm{vol}(U)+\mathrm{vol}(\overline{U})}.\qe
Since this holds for all partitions $U,\overline{U}$, we get
\[h\leq \overline{h}\frac{\mathrm{vol}(U)+\mathrm{vol}(\overline{U}))}
{2\min(\mathrm{vol}(U),\mathrm{vol}(\overline{U}))}=
\overline{h}(\frac{1}{2}+\frac{1}{2\mathcal{R}})\leq
\overline{h}\frac{N}{N-1}\] where we used Lemma \ref{16}. The
remaining inequality follows from Theorem \ref{Lemma2}.
\end{proof}
\begin{coro}\label{coro3} Let $\Gamma$ be an unweighted graph. If there exists a partition
 $U,\overline{U}$ of the vertex set $V$ such that vol($U)=$vol($\overline{U}$)
 then \be \label{coro4} h\leq \overline{h} \leq 1.\qe
\end{coro} If $\Gamma = K_2$, we even have equality in \rf{coro4},
i.e.  $h= \overline{h} = 1$. Note that, in general there does not
exist a partition $U,\overline{U}$ of $V$ such that vol$(U)=$
vol$(\overline{U})$. Counterexamples are regular graphs if $N$ is
odd or so-called wheel graphs $W_N$ with $N$ vertices and degree
sequence $\pi =\{N-1,3,\ldots, 3\}$ if $N-1$ is not a multiple of
$3$.

\begin{figure}\begin{center}
\includegraphics[width =
4.3cm]{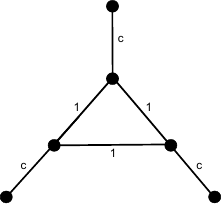}\caption{\label{Fig} For
sufficiently small weights $c$ this graph shows that Lemma
\ref{16} is not true for weighted graphs.}
\end{center}
\end{figure}

\begin{example}\label{23} For a complete graph $K_N$ on $N$ vertices we have
\begin{equation}\label{18} h=  \left\{
\begin{array}{c c} \frac{N}{2(N-1)}& \mbox{$N$ even} \\
 \frac{N+1}{2(N-1)}&\mbox{$N$ odd}
\end{array} \right.\end{equation}
and
\begin{equation}\label{21} \overline{h} = \left\{
\begin{array}{c c} \frac{N}{2(N-1)}& \mbox{$N$ even} \\
 \frac{N+1}{2N}&\mbox{$N$ odd.}
\end{array} \right.\end{equation}
This example shows that, for complete graphs, we have equality in
\rf{11} if $N$ is odd and in \rf{coro4} if $N$ even, respectively.
\end{example}
\noindent\textbf{Weighted graphs}: \\\\ Lemma \ref{16} does not
hold for weighted graphs. This can be seen by considering
sufficiently small weights $c$ in Figure \ref{Fig}.
In particular, it turns out that inequality \rf{29} does not hold
for weighted graphs. However, we have the following result for
weighted graphs.
\begin{lemma}\label{32} Let $\Gamma$ be a weighted graph with $N$
vertices and let $M \leq N$ be defined as in \rf{defM}, then
\begin{equation}
\frac{M-1}{M+1}\leq \mathcal{R}\leq 1.
\end{equation}
\end{lemma}
\begin{proof}
Using the notation from the proof of Lemma \ref{16} we conclude
from \rf{28}:
\begin{eqnarray*}\mathrm{vol}(\overline{U}_N)- \mathrm{vol}(U_M)&=&
\mathrm{vol}(\overline{U}_M)- \mathrm{vol}(U_M) \\&\leq&
\frac{\mathrm{vol}(U_M) +
\mathrm{vol}(\overline{U}_M)}{M}\\&=&\frac{\mathrm{vol}(U_M) +
\mathrm{vol}(\overline{U}_N)}{M}\end{eqnarray*}This implies
\[\frac{M-1}{M+1}\leq \frac{\mathrm{vol}(U_M)}{\mathrm{vol}(\overline{U}_N)}
\leq \frac{\mathrm{vol}(U_N)}{\mathrm{vol}(\overline{U}_N)}\leq
\mathcal{R}.\]
\end{proof}
Note that Figure \ref{Fig} does not contradict Lemma \ref{32} for
all $c>0$. Similarly to Theorem \ref{17} we obtain for weighted
graphs:
\begin{theo}
Let $\Gamma$ be a weighted graph and let $M$ be defined as in
\rf{defM}, then \be \frac{M-1}{M} h\leq
\frac{2\mathcal{R}}{1+\mathcal{R}} h\leq\overline{h}\leq 1.\qe
\end{theo}
\section{Neighborhood graphs}\label{section7}
In the preceding sections we have used geometric properties (like
the Cheeger, the dual Cheeger constant or the diameter) of the
underlying graph in order to control the eigenvalues of the graph
Laplace operator.

In this section, we shall use a new, conceptually different
approach in order to control the eigenvalues of the graph Laplace
operator. Instead of using geometric properties of $\Gamma$
itself, we shall use the geometric properties of the neighborhood
graph $\Gamma[l]$ of $\Gamma$, to be defined shortly, in order to
control the eigenvalues of $\Delta$.

As a motivation, consider  the following result for infinite
graphs \cite{Fujiwara96, Keller10}:
\begin{theo}\label{Fuji}
If $\Gamma$ is a locally finite graph, then \be
\label{Fuji2}1-\sqrt{1-\alpha^2(\Gamma)}\leq \inf
\mathrm{spec}(\Delta) \leq\sup \mathrm{spec}(\Delta)\leq
1+\sqrt{1-\alpha^2(\Gamma)},\qe where \bel{inf1} \alpha(\Gamma) =
\inf_{W\subseteq V,
  |W|<\infty}\frac{E(W,\overline{W})}{\mathrm{vol}(W)}
\qe is a version of the Cheeger constant for an infinite graph.
\end{theo}
Thus, for infinite graphs, it is possible to control the supremum
and the infimum of spec($\Delta$) by the Cheeger constant
$\alpha(\Gamma)$. Clearly, this estimate is not useful for finite
graphs as $\alpha(\Gamma)=0$ in that case, because we may then
simply take $W=V$. The point here is that $\inf
\mathrm{spec}(\Delta)=\lambda_{0} = 0$ for finite graphs, but not
necessarily for infinite graphs, and this is the content of the
lower bound in Theorem \ref{Fuji} in qualitative terms. It is
remarkable that the constant $\alpha(\Gamma)$ at the same time may
also  yield a nontrivial upper spectral bound for
an infinite graph. \\
In the finite graph case we showed in Corollary \ref{62} that a
similar result to Theorem \ref{Fuji} is true if the eigenfunction
that corresponds to the largest eigenvalue is sufficiently
localized. In this section, we shall show (Corollary \ref{coro5})
that it is possible to control the maximal and the smallest
nonzero eigenvalue of a finite graph in a similar way as in
\rf{Fuji2}, if we use the Cheeger constant $h[l]$ of the
neighborhood graph $\Gamma[l]$, for $l$ even, instead of the
Cheeger constant $h$ of the graph $\Gamma$ itself. In particular,
we obtain new lower bounds for the second smallest eigenvalue that
can improve the classical Cheeger estimate \rf{13c}; for a
discussion and examples see the next section.
\begin{defi}
For a graph $\Gamma =(V,E)$ its neighborhood graph $\Gamma[l] =
(V,E[l])$ of order $l\geq1$ is the graph with the same vertex set
$V$ whose  edge set $E[l]$ is defined in the following way: The
weight $w_{ij}[l]$ of the edge $e[l]=(i,j)$ in $\Gamma[l]$ is
given by
\[w_{ij}[l] = \sum_{k_1,\ldots, k_{l-1}}\frac{1}{d_{k_1}}\ldots
\frac{1}{d_{k_{l-1}}}w_{ik_1}w_{k_1k_2}\ldots w_{k_{l-1}j}\] if
$l>1$ and we set $w_{ij}[l]=w_{ij}$ if $l=1$, i.e.
$\Gamma[1]=\Gamma$. In particular, $i$ and $j$ are neighbors in
$\Gamma[l]$ if there exists at least one path of length $l$
between $i$ and $j$ in $\Gamma$.
\end{defi}

\begin{figure}\begin{center}
\includegraphics[width =
3.3cm]{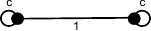}\caption{\label{Fig4}}
\end{center}
\end{figure}

\begin{figure}
\begin{center}
\includegraphics[width =
5.3cm]{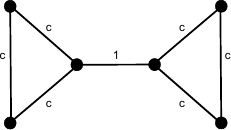}\caption{\label{Fig5}}
\end{center}
\end{figure}

\begin{rem}
\begin{itemize}\item[(i)]The idea of neighborhood graphs is the following: Define a family of
graphs $\Gamma[l]$, $l\geq1$ that encodes the transition
probabilities of the $l$-step random walk on the graph $\Gamma$.
\\
We give an alternative definition of the neighborhood graphs.  The
neighborhood graph $\Gamma[l]=(V, E[l])$ of $\Gamma=(V,E)$ has the
same vertex set V and the weights of the edges of $\Gamma[l]$ are
defined by
$$w_{ij}[l] := P^l(i,j)d_i,$$ where $P^l(i,j)$ is the probability
that a random walker starts at vertex $i$ and moves in $l$ steps
to vertex $j$. \item[(ii)] Neighborhood graphs are directly
related to the discrete heat kernel $p_t(i,j)$ studied in
\cite{Coulhon98}. We have the following relationship:
$$p_t(i,j) = \frac{w_{ij}[t]}{d_id_j}$$
\end{itemize}
\end{rem}
\begin{lemma}
The weights of the neighborhood graphs satisfy the following
semi-group identity:
$$w_{ij}[l] = \sum_t \frac{1}{d_t} w_{it}[k]w_{tj}[l-k]  \text{ for  } 1\leq k<l$$
\end{lemma}
\begin{proof}The proof follows from a direct calculation and we
omit it here.
\end{proof}
In order to become familiar with the concept of neighborhood
graphs we consider the following examples:
\begin{example}\label{examp3}
Consider the family of graphs in Figure \ref{Fig4} for $c\geq 0$.
Let $W$ be the adjacency operator of the graph in Figure
\ref{Fig4}. $W$ can be represented as
\[W=W[1]= \left(\begin{array}{cc} c&1\\1&c
\end{array}\right).\] If we go to higher order neighborhood graphs $\Gamma[l]$ of $\Gamma$,
the topological structure remains the same and the adjacency
operator $W[l]$, $l=2,3,4,5$ can be represented as
\[W[2]= \left(\begin{array}{cc}
 \frac{c^2+1}{1+c}&\frac{2c}{1+c}\\\frac{2c}{1+c}&\frac{c^2+1}{1+c}
\end{array}\right) \qquad W[3]= \left(\begin{array}{cc}\frac{c^3+3c}{(1+c)^2}&
\frac{3c^2+1}{(1+c)^2}\\\frac{3c^2+1}{(1+c)^2}&\frac{c^3+3c}{(1+c)^2}
\end{array}\right)\]
\[ W[4]= \left(\begin{array}{cc}\frac{(c^2+1)^2+4c^2}{(1+c)^3}&
\frac{4c^3+4c}{(1+c)^3}\\\frac{4c^3+4c}{(1+c)^3}&\frac{(c^2+1)^2+4c^2}{(1+c)^3}
\end{array}\right) \qquad W[5]=
\left(\begin{array}{cc}\frac{c (5 + 10 c^2 + c^4)}{(1 + c)^4} &
\frac{1 + 10 c^2 + 5 c^4}{(1 + c)^4}\\ \frac{
 1 + 10 c^2 + 5 c^4}{(1 + c)^4} & \frac{c (5 + 10 c^2 + c^4)}{(1 + c)^4}
\end{array}\right)\]
\end{example}
\begin{example}\label{examp4}
As a second example we consider the family of graphs in Figure
\ref{Fig5}. We have
\[W=W[1]= \left(\begin{array}{cccccc} 0 & c& c& 0& 0& 0 \\c& 0& c& 0& 0& 0\\ c& c& 0& 1& 0& 0 \\ 0& 0& 1& 0& c
& c\\0& 0& 0& c& 0& c \\0& 0& 0& c& c&
  0
\end{array}\right) \text{ and}\] \[W[2]= \left(\begin{array}{cccccc} c/2 + c^2/a& c^2/a& c/2& c/a& 0& 0\\ c^2/a&
 c/2 + c^2/a& c/2& c/a& 0& 0\\c/2& c/2&
 1/a + c& 0& c/a& c/a\\ c/a& c/a& 0& 1/a + c& c/2& c/2\\ 0& 0& c/a& c/2&
 c/2 + c^2/a& c^2/a\\0& 0& c/a& c/2& c^2/a& c/2 +
 c^2/a\end{array}\right),\] where $a:= 1 + 2
 c$.
\end{example}

The the neighborhood graphs $\Gamma[l]$ have the following
property:
\begin{lemma}\label{45}$d_i = d_i[l]$ for all $i\in V$ and $l\geq 1$.
\end{lemma}
\begin{proof}We have \begin{eqnarray*}d_i[l]&=&\sum_jw_{ij}[l]
 = \sum_{k_1,\ldots, k_{l-1}}\frac{1}{d_{k_1}}\ldots
\frac{1}{d_{k_{l-1}}}w_{ik_1}w_{k_1k_2}\ldots
w_{k_{l-2}k_{l-1}}\sum_jw_{k_{l-1}j}\\&=&\sum_{k_1,\ldots,
k_{l-2}}\frac{1}{d_{k_1}}\ldots
\frac{1}{d_{k_{l-2}}}w_{ik_1}w_{k_1k_2}\ldots
\sum_{k_{l-1}}w_{k_{l-2}k_{l-1}}\\&\vdots&\\&=&\sum_{k_1}
w_{ik_1}=d_i.\end{eqnarray*}
\end{proof}
\begin{rem}Lemma \ref{45} implies that
$(\cdot,\cdot)_\mu = (\cdot,\cdot)_{\mu[l]}$ and thus
$\ell^2(V,\mu) = \ell^2(V, \mu[l])$, as will be frequently
utilized below.
\end{rem}
\begin{theo}\label{42} For any function $u\in\ell^2(V,\mu)$ we have
\be \label{43} (I-(I-\Delta)^{l})u= (I-P^{l})u= \Delta[l]u,\qe
where $P$ is the transition probability operator of a random walk
on $\Gamma$ and $\Delta[l]$ is the graph Laplace operator on
$\Gamma[l]$.
\end{theo}
\begin{proof}For any function $u \in \ell^2(V,\mu)$ we have
\begin{eqnarray*}
(I-P^l)u(i) &=&  u(i) - \sum_jq_{ij}u(j),
\end{eqnarray*}
where $q_{ij} := \sum_{k_1,\ldots,
k_{l-1}}\frac{w_{ik_1}}{d_i}\ldots
\frac{w_{k_{l-1}j}}{d_{k_{l-1}}}$ is $ij$-th entry of $P^l$. Thus,
\begin{eqnarray*}
(I-P^l)u(i) &=& u(i) - \sum_j\sum_{k_1,\ldots,
k_{l-1}}\frac{w_{ik_1}}{d_i}\ldots
\frac{w_{k_{l-1}j}}{d_{k_{l-1}}}u(j).
\end{eqnarray*}
Using the definition  $w_{ij}[l] = \sum_{k_1,\ldots,
k_{l-1}}\frac{1}{d_{k_1}}\ldots
\frac{1}{d_{k_{l-1}}}w_{ik_1}w_{k_1k_2}\ldots w_{k_{l-1}j}$ and
$d_i = d_i[l]$, this yields
\begin{eqnarray*}
(I-P^l)u(i) &=&    \frac{1}{d_i[l]}
 \sum_{j}w_{ij}[l] (u(i) -u(j))\\&=&
 \Delta[l]u(i).
\end{eqnarray*}
\end{proof}
\begin{lemma}\label{Lemma6}
Let $\Gamma$ be a graph and $\Gamma[l]$ its neighborhood graph of
order $l$. \begin{itemize} \item[(i)] If $\Gamma$ is connected and
$l$ is even, then $\Gamma[l]$ consists of exactly two connected
components iff $\Gamma$ is bipartite. \item[(ii)] If $l$ is odd,
then $\Gamma[l]$ is bipartite iff $\Gamma$ is bipartite.
Furthermore, $\Gamma[l]$ has the same number of connected
components as $\Gamma$.\item[(iii)] The multiplicity $m_1$ of the
eigenvalue one is an invariant for all neighborhood graphs, i.e.
$m_1(\Delta) = m_1(\Delta[l])$ for all $l\geq 1$.\item[(iv)] The
eigenvalues of $\Delta[l]$ satisfy
\[0= \lambda_0[l]\leq \ldots \lambda_{N-1}[l]\leq
1\] if $l$ is even and \[0= \lambda_0[l]\leq \ldots
\lambda_{N-1}[l]\leq 2\] if $l$ is odd. \item[(v)]If $\lambda \neq
0,2$ then $\lambda[l] = 1-(1-\lambda)^l \rightarrow 1$ as
$l\rightarrow \infty$.
\end{itemize}
\end{lemma}
\begin{proof} $\Gamma$ and $\Gamma[l]$ have the same vertex set, thus both
$\Delta$ and $\Delta[l] = (I-(I-\Delta)^{l})$ have $N=|V|$
eigenvalues. Furthermore, every eigenfunction $u_k$ for $\Delta$
and eigenvalue $\lambda_k$ is also an eigenfunction for
$\Delta[l]$ and eigenvalue $1-(1-\lambda_k)^l$. \\ $(i)$ If $l$ is
even, $\lambda_k[l]=0$ iff $\lambda_k =0$ or $\lambda_k =2$.
Recall that the multiplicity of the eigenvalue zero is equal to
the number of connected components of a graph and $2$ is an
eigenvalue iff the graph is bipartite. Now let $\Gamma$ be a
connected, bipartite graph and let $l$ be even. Then $\lambda[l]
=0$ is twice in the spectrum of $\Delta[l]$ since $\lambda_0=0$
and $\lambda_{N-1}=2$ are in the spectrum of $\Delta$.
Consequently, $\Gamma[l]$ consists of exactly two connected
components. On the other hand if $\Gamma[l]$ consists of exactly
two connected components, $\lambda[l] =0$ is twice in the spectrum
of $\Delta[l]$, and we know that either the eigenvalue $\lambda=0$
is twice in the spectrum of $\Delta$ or $\lambda=0$ and
$\lambda=2$ are both in the spectrum of $\Delta$. Since we assume
that $\Gamma$ is connected, $\lambda=0$ is a simple eigenvalue,
and thus we can conclude  that $2\in \mathrm{spec}(\Delta)$ and
$\Gamma$ is bipartite. \\ $(ii) -(v)$ follow from simple
calculations.
\end{proof}
\begin{rem}By Lemma \ref{Lemma6} (iv) all eigenvalues of $\Delta[l]$ ($l$
even) are less or equal to $1$. From \rf{9} we observe that this
is only possible because $\Gamma[l]$ ($l$ even) contains many (in
fact $N$) loops. In contrast, for a graphs without loops we have
$1< \frac{N}{N-1}\leq \lambda_{N-1}$.
\end{rem}
We now use Theorem \ref{42} to control the eigenvalues of $\Delta$
on $\Gamma$ in terms of geometric properties of its neighborhood
graphs $\Gamma[l]$.
\begin{theo}\label{theo} Let $\mathcal{A}[l]$
be a lower bound for the eigenvalue $\lambda_1[l]$ of $\Delta[l]$,
i.e. $\mathcal{A}[l] \leq \lambda_1[l]$. Then, the eigenvalues of
$\Delta$ satisfy \be \label{39}
1-\left(1-\mathcal{A}[l]\right)^{\frac{1}{l}}\leq \lambda_1 \leq
\ldots\leq \lambda_{N-1}\leq
1+\left(1-\mathcal{A}[l]\right)^{\frac{1}{l}}, \qe if $l$ is even
and \be
\label{44}1-\left(1-\mathcal{A}[l]\right)^{\frac{1}{l}}\leq
\lambda_1 \qe if $l$ is odd.
\end{theo}
\begin{proof}
Let $u_k$, $k\neq0$ be an eigenfunction for $\Delta$. Using
\rf{43} we obtain
\begin{eqnarray*}
1-(1-\lambda_k)^l &=& \frac{(u_k, (I-(I-\Delta)^l)
u_k)_\mu}{(u_k,u_k)_\mu}\\&=& \frac{(u_k,\Delta[l]
u_k)_{\mu[l]}}{(u_k,u_k)_{\mu[l]}}
\\&\geq& \inf_{u: (u,\mathbf{1})_{\mu[l]}=0}\frac{(u, \Delta[l]
u)_{\mu[l]}}{(u,u)_{\mu[l]}}
\\&=&
\lambda_1[l]\geq \mathcal{A}[l],
\end{eqnarray*}
since $(u,\mathbf{1})_\mu =0$ iff $(u,\mathbf{1})_{\mu[l]} =0$.
Alternatively,  since the eigenvalues of $\Delta[l]$ are given by
$1-(1-\lambda_k)^l$ if $\lambda_k$ are the eigenvalues of $\Delta$
we have
\[1-(1-\lambda_k)^{l}\geq \min_{k\neq 0} 1-(1-\lambda_k)^{l}
= \lambda_1[l]\geq \mathcal{A}[l].\] This implies that for all
$k\neq 0$
\[|1-\lambda_k|\leq (1-\mathcal{A}[l])^{\frac{1}{l}}\] if $l$ is
even and
\[1-\lambda_k\leq (1-\mathcal{A}[l])^{\frac{1}{l}}\] if $l$ is
odd.
\end{proof}
As a concrete example, we use the Cheeger inequality \rf{13c} as a
lower bound for $\lambda_1[l]$, i.e. $\mathcal{A}[l] =
1-\sqrt{1-h^2[l]} \leq \lambda_1[l]$, where $h[l]$ is the Cheeger
constant of the neighborhood graph $\Gamma[l]$.
\begin{coro}\label{coro5}The eigenvalues of $\Delta$ on $\Gamma$
satisfy: \be \label{48}1-\left(1-h^2[l]\right)^{\frac{1}{2l}}\leq
\lambda_1\leq \ldots\leq \lambda_{N-1}\leq
1+\left(1-h^2[l]\right)^{\frac{1}{2l}} \qe if $l$ is even and \be
\label{49} 1-\left(1-h^2[l]\right)^{\frac{1}{2l}}\leq \lambda_1\qe
if $l$ is odd.
\end{coro}
We point out that this result is similar to Theorem \ref{Fuji} for
locally finite graphs. However, the difference is that we have to
use here the Cheeger constant of the neighborhood graph instead of
the Cheeger constant itself.

As another  example, we use the well-known estimate $\lambda_1
\geq \frac{\min_{i,j}w_{ij}}{D\mathrm{vol}(V)}$ for the smallest
nontrivial eigenvalue in terms of the diameter and the volume of a
graph, see for instance \cite{Chung97}. This yields the following
estimates:
\begin{coro}\label{73}
All eigenvalues of $\Delta$ satisfy
\[1 -\left(1-
\frac{\min_{i,j}w_{ij}[l]}{D[l]\mathrm{vol}(V[l])}\right)^{\frac{1}{l}}
\leq \lambda_1\leq \ldots \leq \lambda_{N-1} \leq 1 + \left(1-
\frac{\min_{i,j}w_{ij}[l]}{D[l]\mathrm{vol}(V[l])}\right)^{\frac{1}{l}}
\] if $l$ is even and
\[1 -\left(1-
\frac{\min_{i,j}w_{ij}[l]}{D[l]\mathrm{vol}(V[l])}\right)^{\frac{1}{l}}
\leq \lambda_1\]if $l$ is odd.
\end{coro}

\begin{rem}
Some graph properties, like  the discrepancy or the expansion
property of a graph, can be controlled by the quantity $\rho =
\max_{k\neq0}|1-\lambda_k|$ \cite{Chung97}. Corollary \ref{coro5}
and Corollary \ref{73} or more generally Theorem \ref{theo} can be
used to derive explicit bounds for those quantities. As one
particular application, we show in Section \ref{section8} how
Corollary \ref{coro5} can be used to control the convergence of
random walks on graphs.
\end{rem}
We can use Theorem \ref{42} to obtain further eigenvalue
estimates.
\begin{theo}\label{theo3}Let $\mathcal{B}[l]$ be any upper bound for
$\lambda_1[l]$, i.e. $\lambda_1[l] \leq \mathcal{B}[l]$.
 Then the eigenvalues of $\Delta$ satisfy \be \label{cond3}\lambda_1\leq
1-\left(1-\mathcal{B}[l]\right)^{\frac{1}{l}}\qe or
\be\label{cond4} \lambda_{N-1}\geq
1+\left(1-\mathcal{B}[l]\right)^{\frac{1}{l}}\qe if $l$ is even
and \[\lambda_1\leq
1-\left(1-\mathcal{B}[l]\right)^{\frac{1}{l}}\] if $l$ is odd.
\end{theo}
\begin{proof}First note that, by Lemma \ref{Lemma6}, \rf{cond3} and \rf{cond4} are well
defined, if $l$ is even, since we can assume w.l.o.g. that
$\mathcal{B}[l]\leq 1$. Using \rf{43} we obtain $\mathcal{B}[l]
\geq \lambda_1[l] =\min_{k\neq0}1-(1-\lambda_k)^l$. This implies
that for at least one eigenvalue $\lambda_i$, $i\neq 0$ we have
\[(1-\mathcal{B}[l])^{\frac{1}{l}}\leq |1-\lambda_i|\] if $l$ is even and
\[(1-\mathcal{B}[l])^{\frac{1}{l}}\leq 1-\lambda_i\] if $l$ is
odd.
\end{proof}
Using the Cheeger inequality \rf{13a} for $\Gamma[l]$ we obtain:
\begin{coro}\label{coro6} If $l$ is even, and $2h[l]\leq 1$, then we have
\be \label{50} \lambda_1\leq 1-\left(1-2h[l]\right)^{\frac{1}{l}}.
\qe or \be \label{51} \lambda_{N-1}\geq
1+\left(1-2h[l]\right)^{\frac{1}{l}}.\qe If $l$ is odd, we have
\be \label{542}\lambda_1\leq
1-\left(1-2h[l]\right)^{\frac{1}{l}}.\qe
\end{coro}
In the next section we show that the estimate \rf{50} and \rf{542}
for $\lambda_1$ can improve the Cheeger estimate \rf{13a}.
\begin{theo}
\label{theo4} Let $\mathcal{C}[l]$ be any lower bound for the
largest eigenvalue $\lambda_{N-1}[l]$, i.e. $\mathcal{C}[l]\leq
\lambda_{N-1}[l]$. At least one eigenvalue of $\Delta$ is
contained in the interval \be\label{cond5} \left[1 -
\left(1-\mathcal{C}[l]\right)^{\frac{1}{l}},1 +
\left(1-\mathcal{C}[l]\right)^{\frac{1}{l}}\right]\qe if $l$ is
even and the largest eigenvalue of $\Delta$ satisfies
\[\lambda_{N-1}\geq 1 -
\left(1-\mathcal{C}[l]\right)^{\frac{1}{l}}\]if $l$ is odd.
\end{theo}
\begin{proof}Again, by Lemma \ref{Lemma6}, \rf{cond5} is well
defined, if $l$ is even, since $\mathcal{C}[l]\leq\lambda_{N-1}[l]
\leq 1$. Using \rf{43} we have $\mathcal{C}[l]\leq
\lambda_{N-1}[l]=\max_{k} (1-(1-\lambda_k)^l)$. Thus,
\[\min_{k}|1-\lambda_k|\leq
\left(1-\mathcal{C}[l]\right)^{\frac{1}{l}}\]if $l$ is even and
\[\min_{k}(1-\lambda_k)= 1 -\max_k\lambda_k\leq
\left(1-\mathcal{C}[l]\right)^{\frac{1}{l}}\]if $l$ is odd.
\end{proof}
In particular, we have from Theorem \ref{A6} and Theorem
\ref{theo4}:
\begin{coro}\label{coro7} If $l$
is even, and $2\overline{h}[l]\leq 1$,   then at least one
eigenvalue of $\Delta$ is contained in the interval \be \left[1 -
(1-2\overline{h}[l])^{\frac{1}{l}},1 +
(1-2\overline{h}[l])^{\frac{1}{l}}\right].\qe If $l$ is odd, then
the largest eigenvalue satisfies
 \[\lambda_{N-1}\geq 1 -
\left(1-2\overline{h}[l]\right)^{\frac{1}{l}}.\]
\end{coro}

We now turn to the gap phenomenon for eigenvalues, that is, find
some interval that does not contain any eigenvalue.
\begin{theo}\label{theo2}Let $\mathcal{D}[l]$
 be any upper bound for the largest
eigenvalue, i.e. $\lambda_{N-1}[l]\leq \mathcal{D}[l]$. Then all
eigenvalues of $\Delta$ are  contained in the union of intervals
\be\label{cond1}\left[0,1 - \left(1
-\mathcal{D}[l]\right)^{\frac{1}{l}}\right]\bigcup \left[1
+\left(1 -\mathcal{D}[l]\right)^{\frac{1}{l}},2\right]\qe if $l$
is even and the largest eigenvalue of $\Delta$ satisfies
\[\lambda_{N-1}\leq 1 - \left(1
-\mathcal{D}[l]\right)^{\frac{1}{l}}\] if $l$ is odd.
\end{theo}
The proof is similar to the proofs above so we omit it here. We
only note that \rf{cond1} is well defined if $l$ is even since, by
Lemma \ref{Lemma6}, we can assume w.l.o.g. that
$\mathcal{D}[l]\leq1$.

In other words: Let $l$ be even, then for any upper bound
$\mathcal{D}[l]\leq 1$ none of the eigenvalues of $\Delta$ is
contained in the interval $\left(1 - \left(1
-\mathcal{D}[l]\right)^{\frac{1}{l}},1 + \left(1
-\mathcal{D}[l]\right)^{\frac{1}{l}}\right)$. Thus, if $l$ is
even, an upper bound $\mathcal{D}[l]\leq 1$ for $\lambda_{N-1}[l]$
of $\Delta[l]$ can be used to bound all eigenvalues of $\Delta$ on
$\Gamma$ away from $1$. In particular, if $\mathcal{D}[l]<1$,
$\Delta$ then does not possess the eigenvalue $1$, see also Lemma
\ref{Lemma6}. In \cite{BJ}, it was observed that the eigenvalue
$1$ occurs in an unweighted graph whenever there are two nodes
$i_1, i_2$ that are not neighbors themselves, but who possess the
same neighbors, that is, for any $k$, we have $k\sim i_1$ iff
$k\sim i_2$. Thus, a graph satisfying the assumptions of Theorem
\ref{theo2} cannot have any such pair of nodes if
$\mathcal{D}[l]<1$ and $l$ is even.

From the dual Cheeger inequality  \rf{A2} we only obtain an upper
bound for $\lambda_{N-1}$ if $l$ is odd since $\mathcal{D}[l] =
1+\sqrt{1-(1-\overline{h}[l])^2}\geq 1$ for all $l$.
\begin{coro}\label{59} If $l$ is odd, then
\[\lambda_{N-1} \leq 1 + (1-(1-\overline{h}[l])^2)^{\frac{1}{2l}}.\]
\end{coro}
\section[Comparison of the Cheeger estimate with the  estimates] {Comparison of the Cheeger estimates with the estimates
  obtained by the neighborhood graph method}
In the following we compare the Cheeger estimates \rf{13a} and
\rf{13c} with the new estimates  in Corollary \ref{coro5} and
Corollary \ref{coro6}. Recall that Corollary \ref{coro5} and
Corollary \ref{coro6} were obtained by applying \rf{13a} and
\rf{13c} to the neighborhood graph $\Gamma[l]$ and then using the
relationship between the spectrum of $\Delta$ and $\Delta[l]$,
which was established in Theorem \ref{42}.

Comparing Corollary \ref{coro5} with \rf{13c} reveals that our new
estimates improve the Cheeger estimate \rf{13c} if \be \label{46}
h[l]\geq \sqrt{1-(1-h^2)^l},\qe for some $l\geq2$. In general it
is not clear for which graphs $\Gamma$ and which $l$ the equation
\rf{46} is satisfied. However, we can develop some qualitative
intuition about \rf{46}.

We have to distinguish whether $l$ is even or odd. Assume for the
moment that $l$ is even. Clearly, \rf{46} is not satisfied
whenever $\Gamma$ is bipartite since then, by Lemma \ref{45},
$\Gamma[l]$ is disconnected and so $h[l]=0$. In fact,  Corollary
\ref{coro5} yields only the trivial estimate $0\leq \lambda_1$ for
bipartite graphs. In contrast, for graphs that  are not bipartite
the estimate in \rf{48} always yields a non-trivial lower bound
$0<1-\left(1-h^2[l]\right)^{\frac{1}{2l}}\leq \lambda_1$ for the
second smallest eigenvalue $\lambda_1$. A  necessary condition for
strict inequality in \rf{46} is that $h[l]>h$. In order to
understand for which graphs it is likely that this necessary
condition is satisfied we distinguish the following two cases. If
$\lambda_1\geq 2-\lambda_{N-1}$ then $\lambda_1[l]=
1-(1-\lambda_{N-1})^l$. In this case  it is possible that
$\lambda_1[l] <\lambda_1$. Hence, in general,  we cannot expect
that $h[l]>h$ is satisfied, unless the Cheeger estimate on the
graph $\Gamma[l]$ is sharper (this will be made more precise in
the next proposition) than the Cheeger estimate on the graph
$\Gamma$. In particular, if $\lambda_1\geq 2-\lambda_{N-1}$,  we
cannot expect that \rf{46} is satisfied.  On the other hand if
$\lambda_1< 2-\lambda_{N-1}$ then $\lambda_1[l]=
1-(1-\lambda_{1})^l>\lambda_1$ and so it is likely that the
necessary condition $h[l]>h$ is satisfied. Roughly speaking, if
$l$ is even, we can expect that our eigenvalue estimates improve
the Cheeger estimates if the graph in question is closer to
disconnected graph than to bipartite graph.

If $l$ is odd then Corollary \ref{coro5} always yields non-trivial
estimates and $\lambda_1[l]> \lambda_1$ is always satisfied. Thus
we can expect that the necessary condition $h[l]>h$ is, in
general, satisfied if $l$ is odd.

After all, we are mainly interested in the question when the
estimates in Corollary \ref{coro5} improve the Cheeger estimate
\rf{13c}
\begin{satz}\label{53} Let $S[l]:= \frac{1-\sqrt{1-h^2[l]}}{\lambda_1[l]}$ be
the sharpness of the Cheeger estimate \rf{13c} on the graph
$\Gamma[l]$, i.e. the closer $S[l]$ is to $1$ the sharper is the
Cheeger estimate. We set $S[1]=:S$ and $h[1]=:h$. If one of the
following two conditions is satisfied
\begin{enumerate} \item[(i)] $l$ is odd and $S[l]$ satisfies
\be \label{52} S[l] \geq
\frac{1-(1-h^2)^{\frac{l}{2}}}{1-[1-\frac{1}{S}(1-\sqrt{1-h^2})]^l},\qe
 \item[(ii)]  $l$ is even, $S[l]$ satisfies \rf{52}, and
$\lambda_1\leq \lambda_{N-1}$,  \end{enumerate} then the estimates
in Corollary \ref{coro5} improve the Cheeger estimate \rf{13c}.
\end{satz}
\begin{proof}
We have \begin{eqnarray*} S[l] &=&
\frac{1-\sqrt{1-h^2[l]}}{\lambda_1[l]} =
\frac{1-\sqrt{1-h^2[l]}}{1-(1-\lambda_{1})^l}\\ &=&
\frac{1-\sqrt{1-h^2[l]}}{1-[1-\frac{1}{S}(1-\sqrt{1-h^2})]^l}.
\end{eqnarray*} A comparison with \rf{52} yields
\[1-\sqrt{1-h^2[l]}\geq 1-(1-h^2)^{\frac{l}{2}},\] which implies \rf{46}
\end{proof}

\begin{figure}\begin{center}
\includegraphics[width =
7.3cm]{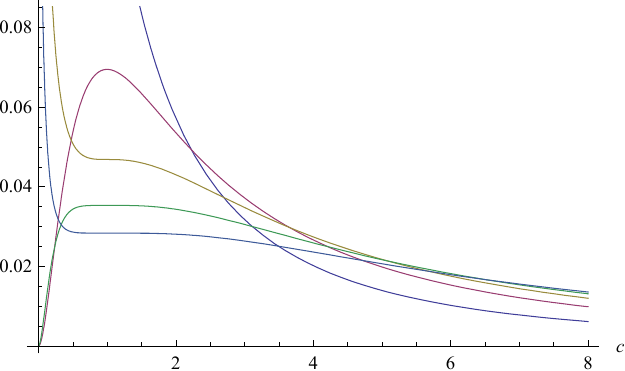}\caption{\label{Fig6} Plot of different lower
bounds $1-(1-h^2[l])^{\frac{1}{2l}}$ for the second smallest
eigenvalue $\lambda_1$ of the family of graphs in Figure
\ref{Fig4}. The Cheeger estimate ($l=1$) is plotted in dark blue,
and the red, yellow, green and light blue curves correspond to
$l=2,3,$4, and $l=5$ respectively.   }
\end{center}
\end{figure}
\begin{figure}\begin{center}
\includegraphics[width =
7.3cm]{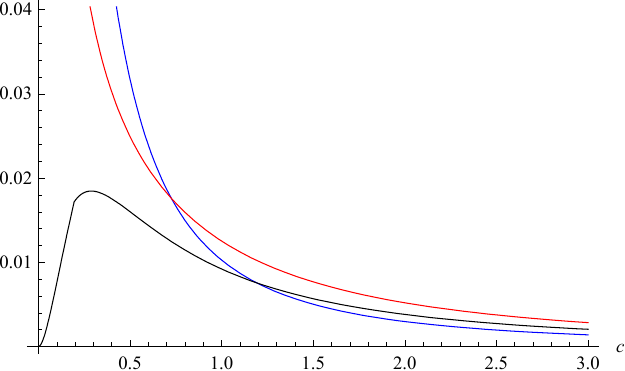}\caption{\label{Fig7} Plot of different lower
bounds $1-(1-h^2[l])^{\frac{1}{2l}}$, for the second smallest
eigenvalue of the family of graphs in Figure \ref{Fig5}. The
Cheeger estimate ($l=1$) is plotted in blue, and the black, red
curves correspond to $l=2$, and $l=3$ respectively.}
\end{center}
\end{figure}\begin{figure}\begin{center}
\includegraphics[width =
7.3cm]{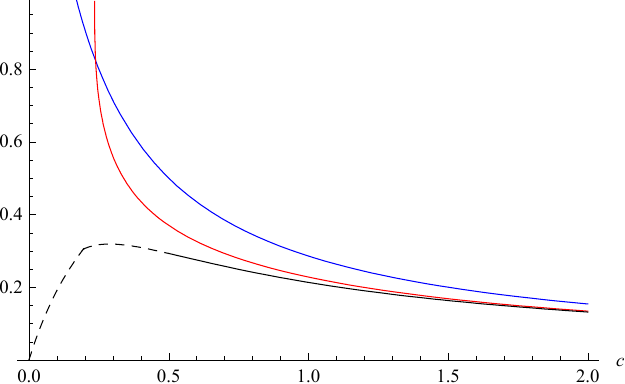}\caption{\label{Fig8}Plot of different upper
bounds $1-(1-2h[l])^{\frac{1}{l}}$, for the second smallest
eigenvalue of the family of graphs in Figure \ref{Fig5}. The
Cheeger estimate ($l=1$) is plotted in blue, and the black, red
curves correspond to $l=2$, and $l=3$ respectively. The dashed
black line indicates that for $c\leq 0.5$ this is not an upper
bound for $\lambda_1$ because in this case \rf{51} holds and
\rf{50} is not satisfied. }
\end{center}
\end{figure}
As an example we consider again the family of graphs in Example
\ref{examp3}. The corresponding Cheeger constants of the
neighborhood graphs are given by $h[1] = \frac{1}{1+c}$, $h[2]=
\frac{2c}{(1+c)^2}$, $h[3]= \frac{3c^2+1}{(1+c)^3}$, $h[4]=
\frac{4c^3+4c}{(1+c)^4}$, and $h[5]= \frac{1+10c^2+c^4}{(1+c)^5}$.
In Figure \ref{Fig6} we plot the lower bounds
$1-(1-h[l]^2)^{\frac{1}{2l}}$, for the second smallest eigenvalue
$\lambda_1$, for different values of $l$. We observe that for
$c<2.2$ the Cheeger inequality ($l=1$) yields the best estimate
for the second smallest eigenvalue $\lambda_1$ of the graph in
Figure \ref{Fig4}. However, if we increase $c$ the estimates for
larger values of $l$ become better. This confirms the intuition
that graphs that are closer to disconnected than to bipartite
graphs are likely to satisfy \rf{46}.

As a second example, we consider the family of graphs in Example
\ref{examp4}. The Cheeger constants of the neighborhood graphs are
given by \linebreak  $h[1] = \min\{\frac{1}{6c+1}, \frac{1}{2}\}$,
$h[2] = \min\{\frac{4c}{(6c+1)(2c+1)}, \frac{5c + 8c^2}{2(1 + 6c +
8c^2)}\}$, and \linebreak $h[3]=
\min\{\frac{12c^2+4c+1}{(1+2c)^2(1+6c)},
\frac{c(12c^2+12c+7)}{(1+2c)^2\min\{8c,2+4c\}}\}$. In Figure
\ref{Fig7} the lower bounds $1-(1-h[l]^2)^{\frac{1}{2l}}$, for the
second smallest eigenvalue $\lambda_1$, are plotted for $l =
1,2,3$. Again, for small $c$ the Cheeger estimate, $l=1$ yields
the best estimate. However, if $c>0.8$ then the estimate for $l=3$
improves the Cheeger estimate.

We can also compare the upper bound for $\lambda_1$ in Corollary
\ref{coro6} with the Cheeger inequality \rf{13a}. We observe that
if $l$ is odd and \be\label{47}h[l]\leq \frac{1-(1-2h)^l}{2}\qe is
satisfied then the estimates in Corollary \ref{coro6} improve the
Cheeger estimate \rf{13a}. If $l$ is even, $h[l]\leq\frac{1}{2}$,
\rf{47} is satisfied,  and  $\lambda_1\leq 2-\lambda_{N-1}$, then
Corollary \ref{coro6} improves the Cheeger estimate \rf{13a}. If
$l$ is even, we have to assume that $\lambda_1\leq
2-\lambda_{N-1}$, because otherwise \rf{51} holds instead of
\rf{50} and so we do not always have an upper bound for
$\lambda_1$. Similarly to Proposition \ref{53} we obtain:
\begin{satz} Let $s[l]:= \frac{\lambda_1[l]}{2h[l]}$ be the
sharpness of the upper Cheeger estimate for the second smallest
eigenvalue $\lambda_1[l]$ of $\Gamma[l]$. If one of the following
two conditions is satisfied

\begin{enumerate}\item[(i)] $l$ is odd and \be \label{54} s[l]\geq \frac{1-(1-s2h)^l}{1-(1-2h)^l}\qe
\item[(ii)] $l$ is even and in addition to \rf{54},
$h[l]\leq\frac{1}{2}$, and $\lambda_1\leq
2-\lambda_{N-1}$,\end{enumerate} then the estimates in Corollary
\ref{coro6} improve the Cheeger estimate \rf{13a}.
\end{satz}
The proof is straightforward so we omit it here.

It turns out that, if we consider the family of graphs in Example
\ref{examp3}, the different upper bounds $\lambda_1\leq
1-(1-2h[l])^{\frac{1}{l}}$ for the second smallest eigenvalue are
the same for all $l$. In contrast for the family of graphs in
Example \ref{examp4} the plot in Figure \ref{Fig8} shows that the
estimate in Corollary \ref{coro6} can improve the Cheeger estimate
\rf{13a}. For example if $c>0.3$ the estimate for $l=3$ improves
the Cheeger estimate \rf{13a}. Comparing Figure \ref{Fig7} and
Figure \ref{Fig8} shows that the estimates in Corollary
\ref{coro5} and Corollary \ref{coro6} improve both Cheeger
estimates \rf{13a} and \rf{13c} at the same time if $c>0.8$ and
$l=3$.

Finally we observe that if $h[l]$ is not contained in the interval
$\left[\frac{1-(1-2h)^l}{2},\sqrt{1-(1-h^2)^l}\right]$ (where
 the interval is the empty set if
$\frac{1-(1-2h)^l}{2}>\sqrt{1-(1-h^2)^l}$) then at least one of
the Cheeger estimates \rf{13a} and \rf{13c} is improved by the
estimates in Corollary \ref{coro5} and Corollary \ref{coro6}.

\begin{figure}\begin{center}
\includegraphics[width =
7.3cm]{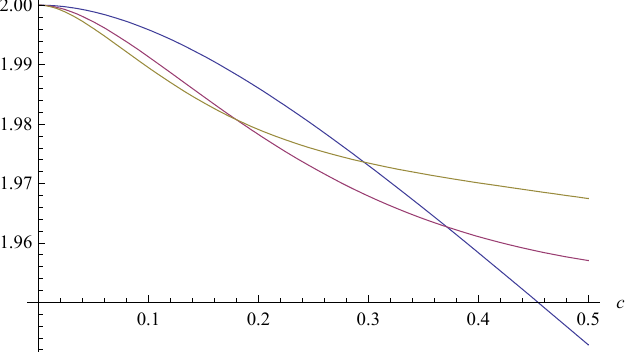}\caption{\label{Fig9} Plot of different
upper bounds $1 + (1-(1-\overline{h}[l])^2)^{\frac{1}{2l}}$ in
Corollary \ref{59} for the largest eigenvalue of the family of
graphs in Figure \ref{Fig4}. The dual Cheeger estimate (Theorem
\ref{A6}) ($l=1$) is plotted in dark blue, and the red and yellow
curves correspond to $l=3$ and $l=5$ respectively. Note that for
the family of graphs in Figure \ref{Fig4} we have $\overline{h}[l]
= h[l]$.}
\end{center}
\end{figure}

Similarly, one can also compare the dual Cheeger estimate for the
largest eigenvalue in Theorem \ref{A6} with the estimates obtained
in Corollary \ref{59}. We observe that, similarly to the Cheeger
estimate, the neighborhood graph method can improve the estimates
obtained from the dual Cheeger estimate. As one such example we
consider the family of graphs in Example \ref{examp3}. In Figure
\ref{Fig9} we plot different upper bounds $1 +
(1-(1-\overline{h}[l])^2)^{\frac{1}{2l}}$ for the largest
eigenvalue $\lambda_{N-1}$. For the largest eigenvalue one can
derive similar results as in \rf{46} and Proposition \ref{53}.
However, we do not want to go into further detail here because the
calculations are exactly the same as before.

We conclude this section by noting that the neighborhood graph
method is very powerful because it may improve any known
eigenvalue estimate.

\section{An example}\label{section}
As discussed above, the highest eigenvalue $\lambda_{N-1}$ of
$\Delta$ becomes largest for bipartite and smallest for complete
graphs, respectively. And a guiding question for this paper is
what can we say about the highest eigenvalue of graphs that are
neither bipartite nor complete, i.e., what structural properties
of $\Gamma$ lead to a highest eigenvalue $\lambda_{N-1}$ close to
2, or very different from 2, respectively.

In order to develop some further intuition about the highest
eigenvalue, we now consider the following example. Let $\Gamma_0$
be a bipartite graph with $N$ vertices. We consider a highest
eigenfunction $\bar{u}$ that is $+1$ on one class and $-1$ on the
other class of vertices, as described above. In particular by
\rf{8k}, \bel{highest3}
\frac{\frac{1}{2}\sum_k\sum_{j}w_{kj}(\bar{u}(j)-\bar{u}(k))^2}{\sum_i
  d_i\bar{u}(i)^2}=2.
\qe By adding another vertex $i_0$ and connecting it to one of the
vertices $i_1$ of $\Gamma_0$ we obtain a new bipartite graph
$\Gamma_1$. We extend $\bar{u}$ by $\bar{u}(i_0)=0$ to $\Gamma_1$.
Thus, the numerator and the denominator of \rf{highest3} are both
increased by $w_{i_0i_1}$. Let $\Gamma_0$ be sufficiently large,
i.e., $\sum_i d_i$ is sufficiently large, then we can achieve for
$\Gamma_1$ that for any given small $\epsilon
>0$, \bel{highest4} \frac{\frac{1}{2}\sum_k\sum_{j}w_{kj}(\bar{u}(j)-\bar{u}(k))^2}{\sum_i
  d_i\bar{u}(i)^2}>2-\epsilon.\qe Now,
this is not affected when we construct a graph $\Gamma$ by
attaching another graph $\Gamma_2$ at $i_0$ and extend $\bar{u}$
by 0 to all of $\Gamma_2$. For instance, $\Gamma_2$ could be a
complete graph $K_M$ with $M$ vertices, for any $M$. In
particular, the difference $2-\lambda_{N-1}$  which has to be
smaller than $\epsilon$ by \rf{LambdaK}, is not very sensitive to
the shape of $\Gamma_2$. This implies, for instance, that
$2-\lambda_{N-1}$ cannot reflect a global quantity like the
clustering coefficient $C$ of \rf{2a.1} that expresses an averaged
difference from a graph being bipartite. In fact, our construction
of attaching a complete graph $K_M$ to a bipartite graph
$\Gamma_0$ through a connecting node produces a graph with $C$
arbitrarily close to its maximal value 1 when $M$ is sufficiently
large.

By extending this example, we can also see that we should have
many eigenvalues $\lambda$ for which $2-\lambda$ is small when the
graph possesses several relatively large bipartite or almost
bipartite parts that are only loosely connected with the rest. (By
\rf{43}, the neighborhood graph $\Gamma[2]$ of such a graph
contains several large components that are only loosely connected,
i.e. many eigenvalues $\lambda_k[2]$ that are small.) This is
analogous to the fact that a graph possesses several small
eigenvalues when it has many relatively large components that are
only loosely connected to the rest, that is, when the graph can be
easily decomposed into several large clusters. Of course, for a
nonconnected graph, that is, one with several components without
links between them, the spectrum simply is the union of the
spectra of the components. Therefore, by the continuity principle,
a graph consisting of clusters that are only loosely connected to
each other has its spectrum approximated by the spectra of these
clusters, that is, by the one of the graph resulting from deleting
the few links between the clusters.

\section{Controlling the largest eigenvalue in terms of a local clustering coefficient}

In this section, we shall provide a more technical estimate from
above for the highest eigenvalue $\lambda_{N-1}$. For that
purpose, we shall first derive some general identity, for a
function $u$ on the vertex set $V$ of $\Gamma$.
\begin{lemma} \label{Lemma1}Let $u$ be an eigenfunction
of $\Delta$ for the eigenvalue $\lambda$. Then, the following
identity holds: \bel{highest2} 2-\lambda =
\frac{(\Delta[2]u,u)_\mu}{(\Delta u,u)_\mu}=\frac{\sum_i
\frac{1}{d_i}\sum_{j,k}w_{ij}w_{ik}(u(j)-u(k))^2}{\sum_i
\sum_{j}w_{ij}(u(i)-u(j))^2}\qe
\end{lemma}
\begin{rem}
This identity follows from Theorem \ref{42} for $l=2$ and \rf{C7b}
applied to the neighborhood graph. Here, as an alternative, we
shall provide a direct proof.
\end{rem}
\begin{proof}
\ba \nonumber
& &\sum_i \frac{1}{d_i}\sum_{j,k}w_{ij}w_{ik}(u(j)-u(k))^2\\
\nonumber &=& \sum_i\frac{1}{d_i}
\left(\sum_{j,k}w_{ik}w_{ij}u(j)^2
-2 \sum_{j,k}w_{ij}w_{ik}u(j)u(k) +\sum_{j,k}w_{ij}w_{ik}u(k)^2\right)\\
\nonumber
&=& 2\sum_i \left(\sum_{j}w_{ij}u(j)^2 -\frac{1}{d_i} (\sum_{j}w_{ij}u(j))^2\right)\\
\nonumber &=& 2\sum_i \sum_{j}w_{ij}u(j)^2 -\sum_i 2d_i
\left(\frac{1}{d_i}\sum_{j}w_{ij}u(j)\right)^2. \ea We now observe
that we can replace $u$ by $u-u(i)$ in the first and hence also in
all subsequent lines. This yields \ba \nonumber
& &\sum_i \frac{1}{d_i}\sum_{j,k}w_{ij}w_{ik}(u(j)-u(k))^2\\
\nonumber &=&2\sum_i \sum_{j}w_{ij}(u(j)-u(i))^2 -\sum_i 2d_i
\left(\frac{1}{d_i}\sum_{j}w_{ij}(u(j)-u(i))\right)^2\\
\nonumber &=&2\sum_i \sum_{j}w_{ij}(u(j)-u(i))^2 -\sum_i 2d_i
(\Delta u(i))^2. \ea Since $u$ is an eigenfunction, $\Delta u
=\lambda u$ for some eigenvalue $\lambda$, then, recalling
\rf{8k}, we obtain \bel{highest1} \sum_i
\frac{1}{d_i}\sum_{j,k}w_{ij}w_{ik}(u(j)-u(k))^2=2\lambda
(2-\lambda)\sum_i d_i u(i)^2. \qe Using \rf{8k} again, we can also
reformulate this as \[ 2-\lambda = \frac{\sum_i
\frac{1}{d_i}\sum_{j,k}w_{ij}w_{ik}(u(j)-u(k))^2}{\sum_i
\sum_{j}w_{ij}(u(i)-u(j))^2}.
\]
\end{proof}

We also observe,  by a reasoning similar to the one for Lemma
\ref{Lemma1}:

\begin{lemma}\label{22} Let $u$ be an eigenfunction of $\Delta$ for the eigenvalue $\lambda$. Then, \bel{highest6}
2-\lambda = \frac{2\sum_i
\sum_{k}w_{ik}\left(\frac{1}{d_i}\sum_{j}w_{ij}(u(j)-u(k))\right)^2}{\sum_i
\sum_{j}w_{ij}(u(j)-u(i))^2}. \qe
\end{lemma}
We now employ \rf{highest2} to interpret $2-\lambda_{N-1}$  as
quantifying how much $\Gamma$ is locally different from being
bipartite. Recall that this quantity is 0  iff $\Gamma$ happens to
be bipartite. Note that (\ref{highest6}) can also be used to
estimate the local difference from being bipartite in terms of
$2-\lambda_{N-1}$.

As discussed above, the (global) clustering coefficient $C$ is not
an appropriate measure for the difference $2-\lambda_{N-1}$. For
instance, in Section \ref{section} we constructed a graph whose
largest eigenvalue is close to $2$ although its clustering
coefficient is close to $1$. However, we shall see that it is
possible to control $2-\lambda_{N-1}$ by the following local
clustering measure \bel{highest5} C_0:=\min_{e=(i,j)}
\frac{\alpha_i + \alpha_j}{2},\qe where
\[\alpha_i := \frac{\sum_{j:e=(i,j)\in \bigtriangleup}w_{ij}}{d_i}\] and
$e=(i,k)\in \triangle$ ($i \in \triangle$) denotes that the edge
$e=(i,k)$ (the vertex $i$) is contained in some triangle. Hence,
$\alpha_i$ is the fraction of weights $w_{ij}$ for fixed $i$ that
are contained in some triangle. In particular, if $i\notin
\bigtriangleup$ then $\alpha_i=0$. Again, $C_0=0$ for a bipartite
and $C_0=1$ for a complete graph. Furthermore, we define
\begin{equation}W := \left(\min_{i\in\bigtriangleup}\min_{k:e=(i,k) \in \bigtriangleup}
\sum_{l:(i,k,l)\in\bigtriangleup}\frac{d_i}{d_l}\frac{w_{li}w_{lk}}{w_{ik}}\right)^{1/2}
\end{equation}   and  \be \overline{d} := \max_{i \in \bigtriangleup}
d_i.\qe
\begin{theo}\label{30} The largest eigenvalue $\lambda_{N-1}$ of $\Delta$ can be
controlled from above by \be 2- \frac{1}{2 \overline{d}}\;C_0
\left(\frac{W}{1 + W}\right)^2 =: 2-H \geq \lambda_{N-1}.\qe
\end{theo}\label{}
\begin{proof}First we rewrite \rf{highest2} for the largest eigenvalue $\lambda_{N-1}$ in the following form:
\be \label{highest9}2-\lambda_{N-1} = \frac{ \sum_{e=(i,j)}
\left(\frac{1}{d_i}\sum_{k}\mu_{ij}w_{ik}(u_{N-1}(j)-u_{N-1}(k))^2
+\frac{1}{d_j}\sum_{k}\mu_{ji}w_{jk}(u_{N-1}(i)-u_{N-1}(k))^2\right)}
{\sum_{e=(i,j)}\mu_{ij}(u_{N-1}(i)-u_{N-1}(j))^2+\mu_{ji}(u_{N-1}(j)-u_{N-1}(i))^2}\qe
where again $\mu_{ij}=w_{ij}$ if $i\neq j$ and $\mu_{ij} =
\frac{1}{2}w_{ij}$ if $i=j$.  In order to control
$2-\lambda_{N-1}$ from below we need to match any term
$\mu_{ij}(u_{N-1}(i)-u_{N-1}(j))^2$ in the denominator by some
term in the numerator of comparable magnitude. Since there is
nothing to match if $i=j$ we can use the weights $w_{ij}$ instead
of $\mu_{ij}$ throughout the proof. For ease of notation we will
drop the index ${N-1}$ in the rest of the proof.

For simplicity, we first match the term $w_{ij}(u(i)-u(j))^2$ in
the denominator with
$\frac{1}{d_i}\sum_{k}w_{ij}w_{ik}(u(j)-u(k))^2$ in the numerator.
Because of the symmetry in $i$ and $j$ the second term in the
numerator and denominator can be treated in the same way.

Let $K_1(i)\subset V$ be the set of all neighbors $k$ of $i$ for
which $e=(i,k)\in \triangle$ and \be(u(j)-u(k))^2\geq
\gamma_{ij}^2 (u(i)-u(j))^2\qe is satisfied. The constant
$\gamma_{ij}$ will be used later on to minimize our upper bound
for $\lambda_{N-1}$. Similarly, let $K_2(i)\subset V$ be the set
of all neighbors $k$ of $i$ which satisfy $e=(i,k)\in \triangle$
and \be\label{Case2}(u(j)-u(k))^2 < \gamma_{ij}^2
(u(i)-u(j))^2.\qe Clearly, \[\sum_{k\in K_1(i)\cup K_2(i)}w_{ik} =
\sum_{k:e=(i,k)\in\bigtriangleup}w_{ik} =\alpha_id_i.\] If
$i\notin \bigtriangleup$ then $\alpha_i = 0$ and $K_1(i) = K_2(i)
= \emptyset$. We distinguish the following two cases:
\\ $(i)$
Assume that $\sum_{k\in K_1(i)}w_{ik} \geq \frac{\alpha_id_i}{2}$
is satisfied for vertex $i$. Consequently, there exists a term \be
\frac{1}{d_i}\sum_{k}w_{ij}w_{ik}(u(j)-u(k))^2\geq \frac{\alpha_i
}{2} \gamma_{ij}^2 w_{ij}(u(i)-u(j))^2 \qe in the numerator. Thus
the term $w_{ij}(u(i)-u(j))^2$ in the numerator is matched.
\\ $(ii)$ Now assume that $\sum_{k\in K_2(i)}w_{ik} \geq \frac{\alpha_id_i}{2}$ is satisfied for vertex $i$.
We can not directly match $ w_{ij}(u(i)-u(j))^2$ by using
\rf{Case2} because it could happen that
$\frac{1}{d_i}\sum_{k}w_{ij}w_{ik}(u(j)-u(k))^2 =0$. However,
Equation \rf{Case2} can be used to find a term of comparable size
in the numerator. \rf{Case2} implies that all neighbors $k$ of $i$
in $K_2(i)$ satisfy \be(u(i)-u(k))^2
> (1-\gamma_{ij})^2 (u(i)-u(j))^2.\qe Now the idea is to use the
terms in the numerator several times in order to match
$w_{ij}(u(i)-u(j))^2$. Multiplying the numerator in \rf{highest2}
by $w_{ij}$ yields \begin{eqnarray}
\label{34}w_{ij}\sum_l\frac{1}{d_l}\sum_{m,k}w_{lm}w_{lk}(u(m)-u(k))^2
&=&\nonumber
w_{ij}\sum_l\frac{1}{d_l}\sum_{k}w_{li}w_{lk}(u(i)-u(k))^2
\\&+& w_{ij}\sum_l\frac{1}{d_l}\sum_{k,m\neq
i}w_{lm}w_{lk}(u(m)-u(k))^2.\end{eqnarray}  Now we will only use
the first term on the r.h.s.. The second term can be used to match
other terms in the denominator. The first term on the r.h.s. of
\rf{34} yields:
\begin{eqnarray*}
&&w_{ij}\sum_l\frac{1}{d_l}\sum_{k}w_{li}w_{lk}(u(i)-u(k))^2
\\&\geq& w_{ij}\sum_l\frac{1}{d_l}\sum_{k\in K_2(i)}w_{li}w_{lk}(u(i)-u(k))^2\\&=& w_{ij}\sum_{k\in K_2(i)}
\sum_{l:(i,k,l)\in\bigtriangleup}\frac{1}{d_l}w_{li}w_{lk}(u(i)-u(k))^2\\&>&
(1-\gamma_{ij})^2w_{ij}(u(i)-u(j))^2 \left(\sum_{k\in
K_2(i)}w_{ik}\right)\min_{k\in
K_2(i)}\sum_{l:(i,k,l)\in\bigtriangleup}\frac{1}{d_l}\frac{w_{li}w_{lk}}{w_{ik}}
\\&\geq&
(1-\gamma_{ij})^2w_{ij}(u(i)-u(j))^2\frac{\alpha_i}{2}\underbrace{\min_{k:
e=(i,k)\in \bigtriangleup}
\sum_{l:(i,k,l)\in\bigtriangleup}\frac{d_i}{d_l}\frac{w_{li}w_{lk}}{w_{ik}}}_{:=A(i)}.
\end{eqnarray*}
 Thus $w_{ij}(u(i)-u(j))^2$ is
matched in the numerator. In order to match all other terms of the
form $w_{ip}(u(i)-u(p))^2$ in the denominator for fixed $i$ we
need to use the terms in the numerator at most $\sum_pw_{ip}=d_i$
times. Note that we can use the second term on the r.h.s. of
\rf{34} in order to match other terms of the form
$w_{mp}(u(m)-u(p))^2$ for $m\neq i$. If some vertex $q$ is not
contained in a triangle we have $\alpha_q = 0$ and thus we do not
need to match the terms $w_{qp}(u(q)-u(p))^2$ for fixed $q$. We
conclude that we used the terms in the numerator at most
$\max_{i\in\bigtriangleup} d_i$ times. Because of the symmetry in
$i$ and $j$ the second term in the numerator \rf{highest9} can be
treated in the same way. We obtain the following estimate:
\begin{eqnarray*}2-\lambda_{N-1} &\geq&
\frac{1}{2\overline{d}} \min_{e=(i,j)}\max_{\gamma_{ij}}
\min\left\{a_{ij}, b_{ij}, c_{ij}, d_{ij}\right\},
\end{eqnarray*}
where
\begin{eqnarray*}a_{ij} &:=&
\frac{\alpha_i+\alpha_j}{2}\gamma_{ij}^2,\\b_{ij} &:=&
\frac{\alpha_i}{2}\gamma_{ij}^2
+\frac{\alpha_j}{2}W^2(1-\gamma_{ij})^2,\\c_{ij}&:=&
\frac{\alpha_i}{2}W^2(1-\gamma_{ij})^2
+\frac{\alpha_j}{2}\gamma_{ij}^2,\\d_{ij}
&:=&\frac{\alpha_i}{2}W^2(1-\gamma_{ij})^2+\frac{\alpha_j}{2}W^2(1-\gamma_{ij})^2,
\end{eqnarray*}
and $W^2 := \min_{i\in\bigtriangleup}A(i)$.

Choosing $\gamma_{ij} = \frac{W}{1+W}$ yields\be 2-\lambda_{N-1}
\geq \frac{1}{2\overline{d}}C_0\left(\frac{W}{1+W}\right)^2.\qe
\end{proof}
\begin{rem}Similarly as the local clustering coefficient, the
Olliver-Ricci curvature on a graph \cite{Oll} is related to the
relative abundance of triangles. In \cite{BJL} it is shown that
the largest eigenvalue of $\Delta$ can be controlled from above by
a lower bound for the Olliver-Ricci curvature on graphs.
\end{rem}

With the scheme developed in this section, the control in the
other direction, that is, estimating the largest eigenvalue from
below, does not quite work, because of the following example.
Consider a graph with many cycles of odd length, but all of them
of length at least 5. Here, $C_0(\Gamma)=0$ as there are no
triangles, but $2-\lambda_{N-1}\neq 0$ because the graph is not
bipartite as
bipartite graphs can only have cycles of even length.\\
However, we can control the largest eigenvalue from below in a
different way, as we have seen in Section \ref{cheeg}.

\section{Random walks on graphs and the convergence to
equilibrium}\label{section8} We have seen that neighborhood graphs
are deeply related to random walks on graphs. Hence it is not
surprising that the techniques developed in section \ref{section3}
and \ref{section7} can be applied to random walks on graphs. We
recall the following theorem for the convergence of random walks
on graphs \cite{Grigoryan09}.
\begin{theo} \label{Grigo}For any function $f\in\ell^2(V,\mu)$, set
\[\overline{f} = \frac{1}{\mathrm{vol}(V)}\sum_jd_jf(j).\] Then for
any positive integer $t$, we have
\be\label{61}\|P^tf-\overline{f}\| \leq \rho^t\|f\|,\qe where
$\rho=\max_{k\neq0}|1-\lambda_k|=\max\{|1-\lambda_1|,
|1-\lambda_{N-1}|\}$ is the spectral radius of the transition
probability operator of a random walk  $P$ and
$\|f\|=\sqrt{(f,f)_\mu}$. Consequently, if $\Gamma$ is connected
and not bipartite, then
\[\|P^tf-\overline{f}\|\rightarrow 0\]as $t\rightarrow \infty$,
i.e. $P^tf$ converges to a constant $\overline{f}$ as
$t\rightarrow \infty$.
\end{theo} We define the equilibrium transition probability operator
$\overline{P}: \ell^2(V,\mu) \to \ell^2(V,\mu)$ as
\[\overline{P}u(i) = \frac{1}{\mathrm{vol}(V)}\sum_jd_ju(j).\]
 For all functions
$f\in\ell^2(V,\mu)$ we have, $\overline{P}f=\overline{f}$  and
thus, by \rf{61}, $P^t$ converges to $\overline{P}$ as
$t\rightarrow\infty$ (if $\Gamma$ is not bipartite). As expected,
the equilibrium transition probability for going from $i$ to $j$
only depends on the degree of vertex $j$ (and the volume of the
graph, i.e. the sum of all degrees). In addition, we define the
equilibrium weighted adjacency operator as
$\overline{W}:=D\overline{P}$, where $D: \ell^2(V,\mu) \to
\ell^2(V,\mu)$ is the multiplication operator defined as \[Du(i) =
d_iu(i).\]  The adjacency operator of $\Gamma[l]$ is given by
$W[l] = DP^l$. Hence, if $\Gamma$ is not bipartite,  $\Gamma[l]$
converges to the equilibrium graph $\overline{\Gamma}$ as
$l\rightarrow \infty$ (in the sense that $W[l]$ converges to
$\overline{W}$). The equilibrium graphs of the families of graphs
studied in Example \ref{examp3} and Example \ref{examp4} can be
represented as
\[ \overline{W}= \left(\begin{array}{cc}\frac{1+c}{2}&
\frac{1+c}{2}
\\\frac{1+c}{2}&\frac{1+c}{2}
\end{array}\right)\] and \[\overline{W}= \left(\begin{array}{cccccc}
a_1 & a_1&a_2& a_2&a_1& a_1 \\a_1 & a_1&a_2& a_2&a_1& a_1  \\
a_2 & a_2&a_3& a_3&a_2& a_2 \\a_2 & a_2&a_3& a_3&a_2& a_2\\a_1 &
a_1&a_2& a_2&a_1& a_1 \\a_1 & a_1&a_2& a_2&a_1& a_1
\end{array}\right),\] where $a_1=\frac{4c^2}{12c+2},
a_2=\frac{4c^2+2c}{12c+2}$, and $a_3=\frac{4c^2+4c+1}{12c+2}$.

If $\Gamma$ is bipartite then $\Gamma[l]$ does not converge as
$l\rightarrow\infty$. This can for example be seen from Lemma
\ref{45}. For a bipartite graph $\Gamma$,  $\Gamma[l]$ is
disconnected and not bipartite whenever $l$ is even and
$\Gamma[l]$ is connected and bipartite whenever $l$ is odd.
However, for a bipartite graph $\Gamma$, the subsequence of
neighborhoods graphs $\Gamma[l]$ for $l$ even converges to
$\overline{\Gamma}_{l \text{ even}}$ as $l\rightarrow\infty$. The
corresponding operator $\overline{W}_{l \;even}$ can be
represented as
\[\overline{W}_{l \;even}=
\left(\begin{array}{cc}\overline{W}_1&0\\0&\overline{W}_2\end{array}\right),\]
where $(\overline{W}_k)_{ij}=\frac{2d_id_j}{\mathrm{vol}(V)} =
\frac{d_id_j}{\mathrm{vol}(V_k)}$, for $k=1,2$, if $i$ and $j$
belong to the same subset $V_k$ and $V_1,V_2$ yields a bipartite
decomposition of the vertex set $V$. Thus, $\overline{\Gamma}_{l
\text{ even}}$ is the disjoint union of two complete graphs of
size $|V_1|$ and $|V_2|$. Similarly, if $\Gamma$ is bipartite, the
subsequence $\Gamma[l]$ for $l$ odd converges to
$\overline{\Gamma}_{l \text{ odd}}$ as $l\rightarrow\infty$. In
this case the corresponding operator $\overline{W}_{l\; odd}$ can
be represented as
\[\overline{W}_{l\; odd}=
\left(\begin{array}{cc}0&\overline{W}_1\\\overline{W}_2&0\end{array}\right),\]
where $(\overline{W}_k)_{ij}=\frac{2d_id_j}{\mathrm{vol}(V)}$, for
$k=1,2$, if $i$ and $j$ belong to different subsets. Thus,
$\overline{\Gamma}_{l \text{ odd}}$ is the complete bipartite
graph that has  same bipartite decomposition $V_1,V_2$ of vertex
set $V$ as $\Gamma$.

From Theorem \ref{Grigo} we see that we need to control the
spectral radius $\rho$ of the transition probability operator $P$.
Our results in Section \ref{section3} and \ref{section7} allow us
to control $\rho$ in various different ways. For example,
Corollary \ref{coro5} implies that
\[\rho = \max_{k\neq0}|1-\lambda_k|\leq
(1-h^2[l])^{\frac{1}{2l}}, \text{ for all } l \text{ even.}\] This
gives us the following explicit estimates for Theorem \ref{Grigo}:

\begin{theo}
\[\|P^tf-\overline{f}\| \leq
(1-h^2[l])^{\frac{t}{2l}}\|f\|,\] where we can use the Cheeger
constant $h[l]$ for any even $l$.
\end{theo}
\begin{rem}Instead of considering the convergence in the
norm $\|f|\| = \sqrt{(f,f)_\mu}$, as in \rf{61}, one could also
study stronger notions of convergence, e.g. the relative pointwise
distance \cite{Sinclair93} or other measures of convergence as the
mixing time \cite{Grigoryan09, Sinclair93}. All these quantities
can be bounded from above in terms of the spectral radius of the
transition probability operator $P$. Thus,  the techniques
developed in Section \ref{section3} and \ref{section7} yield
explicit bounds for the convergence of a random walk measured by
any of these notions  of convergence.
\end{rem}

\section{Synchronization in coupled map lattices}
In this section, we present another application of our eigenvalue
estimates.

We consider a coupled map lattice supported by a graph $\Gamma$,
that is, a dynamical system updated at discrete times $t\in
\mathbb{N}$ and of the form \be \label{CML}x_i(t+1) = f(x_i(t)) +
\frac{\epsilon}{d_i}\sum_{j}w_{ij}(f(x_j(t)) - f(x_i(t)))),\qe
where $\epsilon\geq 0$ is the overall coupling strength and
$w_{ij}$ is the strength of the interaction between unit $i$ and
unit $j$. It was discovered by Kaneko \cite{Kaneko90} that the
system \rf{CML} can (asymptotically) synchronize, i.e., $|x_i(t) -
x_j(t)|\rightarrow 0$ for $t\rightarrow \infty$ and all $i,j$,
even if the function $f$ displays chaotic behavior, i.e. its
Lyapunov exponent $\mu(f)$ satisfies \footnote{If $f$ preserves a
``reasonable" good measure, then the Lyapunov exponent does not
depend on the initial value $s(0)$ and a positive Lyapunov
exponent implies positive topological entropy, which is usually
used in the mathematical literature to characterize chaotic
behavior \cite{Robinson}.  In the following we assume that $f$
preserves such a reasonable measure.} \be \mu(f) =
\lim_{T\rightarrow \infty}\frac{1}{T} \sum_{t=0}^{T-1}\ln
|f'(s(t))|>0.\qe Here $s(t)$ is a synchronous solution, i.e.
$x_i(t) = s(t)$ for all $i$.

More precisely, the system \rf{CML} synchronizes under suitable
conditions that will depend  on $\mu(f)$, $\epsilon$, and the
properties of $\Gamma$. In particular, we have the following
criterion for the asymptotic stability of a synchronized state
(which, when fulfilled, implies that \rf{CML} will asymptotically
synchronize when its initial values are sufficiently close to that
state).

\begin{theo}[\cite{Jost01}]\label{sync}
A synchronized state $s(t)$ of the  coupled map lattice \rf{CML}
is asymptotically stable if
\begin{equation}
\label{loc9} {{1-e^{-\mu(f)}} \over {\lambda_1}} < \epsilon <
{{1+e^{-\mu(f)}} \over {\lambda_{N-1}}}.
\end{equation}
Thus, there exists a  range of values of $\epsilon$ for which we
have asymptotic stability if \be \frac{\lambda_{N-1}}{\lambda_1} <
\frac{e^{\mu(f)}+1}{e^{\mu(f)}-1} \mbox{ \qquad and \quad
$\mu(f)>0$}\qe or \be \frac{\lambda_{N-1}}{\lambda_1} >
\frac{e^{\mu(f)}+1}{e^{\mu(f)}-1} \mbox{ \qquad and \quad
$\mu(f)<0$}.\qe The nontrivial case here is, of course, the one
where $\mu(f)>0$.
\end{theo}

Thus, we can determine conditions under which system \rf{CML}
synchronizes.\\
The main point of Theorem \ref{sync} is that the graph in question
should be sufficiently different from both a disconnected graph
(as characterized by $\lambda_1=0$) and a bipartite one (as
characterized by $\lambda_{n-1}=2$). A disconnected graph cannot
synchronize dynamics because its components do not interact.
Dynamics on a bipartite graph need not synchronize because the two
classes can exchange their states every other period, that is, the
bipartite graph can sustain nonsynchronized period 2 oscillations.
For a more general framework for synchronization of coupled
dynamics, see \cite{BAJ}.

In order to apply Theorem \ref{sync}, we need to control the ratio
$\frac{\lambda_{N-1}}{\lambda_1}$. So far, only suitable bounds,
in terms of graph invariants, where known for $\lambda_1$. Using
our estimates for the largest eigenvalue $\lambda_{N-1}$ we can
now control the ratio $\frac{\lambda_{N-1}}{\lambda_1}$ in an
appropriate way. In particular, the our results imply:
\begin{coro}For every graph we have:
\be \frac{\overline{h}}{h}\leq \frac{\lambda_{N-1}}{\lambda_1}\leq
\frac{\min\left[\min_{l\in \N,  \;even}
1+(1-h[l]^2)^{\frac{1}{2l}}, \min_{l\in \N,  \;odd}1 +
(1-(1-\overline{h}[l])^2)^{\frac{1}{2l}}\right]}{\max_{l\in\N}\left(
1-(1-h[l]^2)^{\frac{1}{2l}}\right)}\qe
\end{coro}

\textbf{Acknowledgements.} FB wants to thank Matthias Keller and
Daniel Lenz for stimulating discussions during his stay in Jena.

\end{document}